\documentclass[12pt]{amsart}

\usepackage{amscd}
\usepackage{amssymb}
\usepackage{amsmath}

\newcommand{\C}{{\mathbb{C}}}        
\newcommand{\R}{{\mathbb{R}}}        
\newcommand{\Z}{{\mathbb{Z}}}        
\renewcommand{\H}{{\mathcal{H}}}  
\newcommand{\hyp}{\mathcal{H}}

\newcommand{\G}{{\mathcal{G}}}       

\newcommand{\oddots}{{\mathinner{\mkern1mu\raise1pt\vbox{\kern7pt\hbox{.}}\mkern2mu\raise4pt\hbox{.}\mkern2mu\raise7pt\hbox{.}\mkern1mu}}}

\newcommand{\As}{{(A,\sigma)}}

\newcommand{\Fx}{{F^{\ast}}}
\newcommand{\Fxsq}{{F^{\ast2}}}

\newcommand{\Fsq}{{\Fx/\Fxsq}}

\newcommand{\Lx}{{L^{\ast}}}

\newcommand{\Kx}{{K^{\ast}}}

\newcommand{\qform}[1]{{\langle{#1}\rangle}}                   
\newcommand{\pform}[1]{{\ll\negthinspace{#1}\negthinspace\gg}} 

\newcommand{\IF}[1]{{I^{#1}\negthinspace F}}                   

\newcommand{\IL}[1]{{I^{#1}\negthinspace L}}

\newcommand{\basemu}{\boldsymbol{\mu}}
\newcommand{\mmu}[1]{\basemu_#1}     
\newcommand{\mmut}[1]{{\basemu^{\otimes 2}_#1}}     

\DeclareMathOperator{\Spin}{Spin}           
\newcommand{\Ga}{\mathbb{G}_a}

\newcommand{\semid}{\rtimes}

\newcommand{\Gal}{{\mathrm{{\mathcal{G}}a\ell}\,}}

\DeclareMathOperator{\disc}{disc}

\newcommand{\End}[1]{{\mathrm{End}_{#1}}}
\newcommand{\EndF}{\End{F}}

\newcommand{\Aut}[1]{{{\mathrm {Aut}}_{#1}}}
\DeclareMathOperator{\aut}{Aut}

\newcommand{\Nrd}{{\mathrm{Nrd}}}

\newcommand{\longto}{\longrightarrow}

\newcommand{\injects}{\hookrightarrow}

\newcommand{\ra}{\rightarrow}

\newcommand{\Gm}{\mathbb{G}_m}

\newcommand{\cc}[1]{#1^\circ}

\newcommand{\Asi}{(A_i, \s_i)}

\newcommand{\QZt}{\mathbb{Q}/\mathbb{Z}(2)}
\renewcommand{\mmu}[1]{\basemu_{#1}}
\renewcommand{\mmut}[1]{\basemu_{#1}^{\times 3}}
\newcommand{\mmud}[1]{\basemu_{#1}^{\otimes 2}}

\newcommand{\ach}{\check{\alpha}}

\newcommand{\Dch}{\check{\Delta}}
\newcommand{\dch}{\check{\delta}}

\newcommand{\pr}{\mathbb{P}}

\newcommand{\s}{\sigma}

\renewcommand{\mmu}[1]{\boldsymbol{\mu}_{#1}}

\newcommand{\ua}{\underline{a}}

\newcommand{\ub}{\underline{b}}
\newcommand{\ut}{{\underline{t}}}
\newcommand{\ul}{{\underline{\ell}}}

\newcommand{\oD}{{^1\hskip-.2em D_4}}
\newcommand{\oDe}{{^1\hskip-.2em D_8}}

\newcommand{\oDn}{{^1\hskip-.2em D_n}}
\newcommand{\tDn}{{^2\hskip-.2em D_n}}

\newcommand{\tAn}{{^2\hskip-.2em A_n}}
\newcommand{\tAs}{{^2\hskip-.2em A_6}}

\newcommand{\oEvi}{{^1\hskip-.2em E_6}}
\newcommand{\dEvi}{{^2\hskip-.2em E_6}}

\newcommand{\agrp}[2]{^#1\hskip-.2em #2}

\newcommand{\SN}{\mathrm{SN}}

\newcommand{\e}{\varepsilon}

\newcommand{\tbtmat}[4]{\left( \begin{array}{cc} #1&#2 \\ #3&#4 \end{array} \right) }
\newcommand{\stbtmat}[4]{\left( \begin{smallmatrix} #1&#2 \\ #3&#4 \end{smallmatrix} \right) }

\newcommand{\jordmat}[6]{\left( \begin{array}{ccc} #1&#6&\cdot \\ \cdot&#2&#4 \\ #5&\cdot&#3 \end{array} \right) }
\newcommand{\sjordmat}[6]{\left( \begin{smallmatrix} #1&#6&\cdot \\ \cdot&#2&#4 \\ #5&\cdot&#3 \end{smallmatrix} \right) }

\newcommand{\basjord}{\jordmat{\e_0}{\e_1}{\e_2}{c_0}{c_1}{c_2}}
\newcommand{\sbasjord}{\sjordmat{\e_0}{\e_1}{\e_2}{c_0}{c_1}{c_2}}

\newcommand{\basmat}{\tbtmat{\alpha}{j}{j'}{\beta}}

\newcommand{\trio}[1]{#1^{\times\negthinspace 3}}
\newcommand{\trion}[1]{(#1_0, #1_1, #1_2)}

\newcommand{\diag}{\mathrm{diag}\,}

\newcommand{\Inv}{\mathrm{Inv}\,}

\renewcommand{\aut}{\Aut{}\,}

\newcommand{\Rel}{\mathrm{Rel}\,}

\newcommand{\GOP}{\cc{GO}}

\newcommand{\D}{\Delta}

\renewcommand{\C}{\mathfrak{C}}

\renewcommand{\G}{\Gamma}

\newcommand{\h}{\mathfrak{H}}
\newcommand{\n}{\mathfrak{n}}

\newcommand{\la}{\lambda}

\newcommand{\Cn}{(\C,\n)}

\newtheorem{thm}{Theorem}[section]        
\newtheorem{lem}[thm]{Lemma}
\newtheorem{cor}[thm]{Corollary}
\newtheorem{prop}[thm]{Proposition}

\newtheorem{mainthm}[thm]{Main Theorem}

\newtheorem{ferrlem}[thm]{Ferrar's Lemma}
\newtheorem{movlem}[thm]{Moving Lemma}

\theoremstyle{definition}

\newtheorem{defn}[thm]{Definition}
\newtheorem{eg}[thm]{Example}

\theoremstyle{remark}

\newtheorem{rmk}[thm]{Remark}

\newenvironment{comm}{\noindent{\sl Comments:}}{\hfill$\qed$\medskip}

\newenvironment{pf}{\noindent{\sl Proof:}}{\hfill$\qed$\medskip}
\newenvironment{pf*}{\noindent{\sl Proof:}}{\medskip}

\numberwithin{equation}{section}

\makeatletter
\newenvironment{neqn}%
{\setcounter{equation}{\value{thm}}\begin{eqnarray}}%
{\end{eqnarray}{\stepcounter{thm}}\global\@ignoretrue}
\makeatother

\makeatletter
{\setcounter{equation}{\value{thm}}}%
{{\stepcounter{thm}}\global\@ignoretrue}
\makeatother

\makeatletter
{\setcounter{equation}{\value{thm}}}%
{{\stepcounter{thm}}\global\@ignoretrue}
\makeatother

\makeatletter
\newenvironment{borel}[1]%
{\smallskip \noindent\refstepcounter{thm}{\bf \thethm.}{\bf{ #1.}}}%
{\smallskip \global\@ignoretrue}
\makeatother

\makeatletter
\newenvironment{borel*}%
{\smallskip \noindent\refstepcounter{thm}{\bf \thethm.}}%
{\smallskip \global\@ignoretrue}
\makeatother

\begin{document}

\title[Rost invariant]%
{The Rost invariant has trivial kernel for quasi-split groups of low rank}
\author{R.~Skip Garibaldi}  
\date{30 October 2000}

\begin{abstract}
For $G$ an almost simple simply connected algebraic group defined over
a field $F$, Rost has shown that there exists a canonical map $R_G \!:
H^1(F, G) \ra H^3(F, \QZt)$.  This includes the Arason invariant for
quadratic forms and Rost's mod 3 invariant for Albert algebras as special
cases.  We show that $R_G$ has trivial kernel if $G$ is quasi-split of
type $E_6$ or $E_7$.  A case-by-case analysis shows that it has
trivial kernel whenever $G$ is quasi-split of low rank.
\end{abstract}

\maketitle

\tableofcontents

\bigskip

For $G$ an almost simple simply connected algebraic group over a field $F$, the set of all natural transformations of functors
\[
H^1(?, G) \longto H^3(?, \QZt)
\]
is a finite cyclic group \cite[\S 31]{KMRT} with a canonical
generator.  (Here $H^i(?,M)$ is the Galois cohomology functor which
takes a field extension of your base field $F$ and returns a group if
$M$ is abelian and a pointed set otherwise.  When $F$ has
characteristic zero, $\QZt$ is defined to be $\lim_{\ra} \mmud{n}$ for
$\mmu{n}$ the algebraic groups of $n$th roots of unity; see
\cite[p.~431]{KMRT} or \cite[I.1(b)]{Gille:inv} for a more complete
definition.)  This generator is called 
the {\em Rost invariant} of $G$ 
and we denote it by $R_G$.  In an abuse of notation, we also write
$R_G$ for the map $H^1(F, G) \longto H^1(F, \QZt)$.

This map provides a useful invariant for algebraic structures
classified by $H^1(F, G)$, and an important and typically difficult
question is to describe the kernel of $R_G$.  For example, when $G$ is
split of type $D_n$, $R_G$ is essentially the Arason invariant $\IF3
\ra H^3(F, \mmud{2})$ for
quadratic forms, where $\IF{n}$ is as usual the $n$th power of the
ideal $IF$ of even-dimensional quadratic forms in the Witt ring of $F$.
That
the kernel of the Arason invariant is
precisely $\IF4$ is a quite difficult result due independently to
Merkurjev-Suslin \cite{MS:norm3} and Rost.  (The proof of the main
result of this paper somehow boils  down to this one fact.)
In general, one doesn't even know if the kernel of $R_G$ is trivial.
On the other hand, the question becomes tractable if we 
assume that $G$ is quasi-split.  Generally $R_G$ has nontrivial
kernel; we give easy examples where $G$ is
split of type $D_8$ (in \ref{1D8.eg}) and $B_7$ (in \ref{B7.eg}), and
quasi-split of type $\tAs$ (in \ref{2A6.eg}).
It should be mentioned that $R_G$ can have nontrivial kernel when $G$
is split of type $E_8$ as well; Gille \cite{Gille:E8} has 
produced an example by applying  
his results from \cite{Gille:inv} to reduce the question to the same one for 
a split 
group of type $D_8$.

The principal result in this paper is to enlarge the list of 
quasi-split groups for which the Rost invariant is known to have
trivial kernel. 
\begin{mainthm} \label{MT}
Suppose that $G$ is a quasi-split simply connected group of type 
$E_6$ or $E_7$.  Then
the Rost invariant $R_G$ has trivial kernel.
\end{mainthm}

\begin{borel*}
There are some easy consequences of this theorem that may help
the reader place it in context.  The first is that as a vastly less
powerful corollary, we obtain Serre's ``Conjecture II'' for
quasi-split groups of type $E_6$ and $E_7$, in that if $F$ has
cohomological $p$-dimension $\le 2$ for $p = 2, 3$ (see
\cite[I.3]{SeCG} for a definition), then the main theorem implies that
$H^1(F, G)$ is trivial.  This conjecture appeared in print back in
1962 \cite{Se62}, and remained open for such groups until the 1990s, when
Chernousov (unpublished) and Gille \cite{Gille:sc} proved it (amongst
other cases) independently and by different methods.  Here we get it
for free from the Main Theorem.
\end{borel*}

\begin{borel*}
Another consequence is the following: Suppose that $L$
is a field extension of $F$ of degree relatively prime to 2 and 3 and
that $G$ is a group of type $E_6$ or $E_7$.  Serre asked in
\cite[p.~233, Q.~1]{SeCG:p} if the natural map $H^1(F, G) \ra H^1(L,
G)$ is injective.  Our Main Theorem gives the partial answer that it
has trivial kernel in the case where $G$ is quasi-split.  This result
was already known by experts in the area using arguments special to
groups of type $E_6$ and $E_7$, but as for Conjecture II we get it
for free here.
\end{borel*}

\begin{borel*}
There is also an application to finite-dimensional algebras.  There is
a large family of nonassociative algebras with involution called {\em
structurable 
algebras} which includes central simple associative algebras with
involution (as studied in \cite{KMRT}) and Jordan algebras (with
involution the identity), see \cite{A:survey} for a survey.  The
simple structurable algebras have all been classified, and they
consist (roughly) of the two families already mentioned plus four others.  I
refer to one of these classes, which consists of 56-dimensional
algebras all of which are isomorphic over a separably closed field and
have automorphism group which is simply connected of type $E_6$, as
{\em Brown algebras}.  Now there is a natural equivalence relation
defined on the set of structurable algebras called {\em isotopy}
\cite{AH} which
is weaker than isomorphism, and in the case of Jordan algebras is the
same as the traditional notion of isotopy.  For Albert algebras, it is
known that any algebra 
isotopic to the split one is actually split.  (This is equivalent to
the cohomological statement that the map $H^1(F, F_4) \ra H^1(F, E_6)$
induced by the embedding $F_4 \ra E_6$ described in \ref{E6.rm} has
trivial kernel.)  The Main Theorem here
combined with \cite[4.16(2), 5.12]{G:struct} shows that the same
conclusion holds for Brown algebras, which was previously unknown.
(This has the cohomological interpretation that the map $H^1(F,
E^K_6) \ra H^1(F, E_7)$ induced by the embedding $E^K_6 \ra E_7$
described in \ref{E7.sm} has trivial kernel.)
\end{borel*}

\medskip

The material in \cite{KMRT} is sufficient
to show that the kernel of the Rost invariant
is trivial for quasi-split groups of type $G_2$,
$D_4$ (including those of trialitarian type \cite[40.16]{KMRT}), and
$F_4$, at least away 
from the ``bad primes'' 2 and 3.  As easy corollaries to results
needed for the $E_6$ and $E_7$ cases, we get analogous results for
groups of type $\tAn$, $B_n$, and nontrialitarian groups of type
$D_n$ with small $n$ in Section \ref{BDsec}.
So since $H^1(F, G)$ is always trivial for $G$ split of type $A_n$ or
$C_n$, we get the following:

\begin{thm}
Suppose that $G$ is an almost simple simply connected algebraic group.
If $G$ is
\begin{itemize}
\item quasi-split of (absolute) rank $\le 5$;
\item quasi-split of type $B_6$, $D_6$, or $E_6$; or
\item split of type $D_7$ or $E_7$,
\end{itemize}
then the Rost invariant $R_G$ has trivial kernel.
\end{thm}
\noindent The proofs of these theorems that we will give here and the
material in \cite{KMRT} rely on the ground field having ``good''
characteristic, meaning for our purposes $\ne 2, 3$.  However, it is a
consequence of Gille's main theorem in \cite{Gille:inv} that one
only needs to prove that the Rost invariant has trivial kernel for
fields of characteristic 0.  Consequently, all fields considered here
will be assumed to have characteristic $\ne 2, 3$, but our two theorems
will still hold for all characteristics.  (Of course, in prime
characteristic the group $\QZt$ must be defined somewhat differently
\cite{Gille:inv}, but this affects neither the statemtn of the
theorems nor our proofs.)

Section \ref{BDsec} dispenses with the classical groups.  
(Some of that material is useful later.)
Sections
\ref{foldsec} and \ref{smallsec} contain the material necessary to
reduce questions about the Rost invariant for a larger group to a
subgroup.  That material easily reduces the proof of the main theorem
to considering the quasi-split $\dEvi$ case, which is treated in the
remaining Sections \ref{Gsec} through \ref{E6sec}.\footnote{After this
paper was released as a preprint, Chernousov sent to me a different
proof of the $\dEvi$ case \cite{Ch:rinvlet}, which uses a completely
different argument.  His proof will be published separately.}

\begin{rmk}[Noninjectivity for $F_4$] \label{noninj.F4}
We caution the reader that even when the Rost invariant has trivial
kernel, it may be far from injective.  For example, for $F_4$ the
split group of type $F_4$, the set $H^1(F, F_4)$ classifies 
Albert $F$-algebras.  From known facts about Albert algebras, it is
easy to show that two classes $\alpha_1$, $\alpha_2$ corresponding to
isotopic Jordan algebras $J_1$, $J_2$ have the same Rost invariant.
Since there are many isotopic Albert algebras which are not
isomorphic (for example, over $\R$ there are 3 isomorphism classes of
Albert algebras and two of these are isotopic \cite[p.~119]{Jac:ex}),
the Rost invariant for $F_4$ has trivial kernel but is typically not
injective.   
\end{rmk}

\subsection*{Notations and conventions} 
We say that an algebraic group $G$ is (absolutely) almost simple if it
has finite 
center and no noncentral closed normal subgroups.  When we say that a group
is ``of type $T_n$'', we implicitly mean that it is almost simple of
that type.  We will use the standard notations $\Gm$, $\Ga$, and
$\mmu{n}$ for the algebraic groups with $F$-points $\Fx$, $F$, and the $n$th
roots of unity in $F$, and $\cc{G}$ will always denote the identity
component of an algebraic group $G$.  For a variety $X$ we write $X(F)$ for
its $F$-points.

Our notation for quadratic forms will follow the standard reference
\cite{Lam}.  However, the reader
should be warned of two quirks:
We use the Pfister-approved notation for Pfister forms, so $\pform{a_1,
  \ldots, a_n} := \qform{1, -a_1} \otimes \cdots \otimes \qform{1,
  -a_n}$, and we write $\hyp$ for the hyperbolic plane $\qform{1, -1}$.

The standard reference for Galois cohomology is \cite[\S I.5]{SeCG},
and for algebras with involution (including the groups $\Spin\As$,
$O\As$, and $SO\As$) it is \cite{KMRT}.

\section{Quasi-split groups of type $A$, $B$, and $D$} \label{BDsec}

As indicated in the introduction, the Rost invariant should have trivial kernel
for quasi-split groups of small rank.  To prove this for $E_6$, we
will need a result on groups of type $D$, which also easily settles
this question for groups of type $A$ and $B$.  
(For the results in this
section, our global hypothesis that our fields have characteristic
$\ne 3$ is not required; we need only assume characteristic $\ne
2$.)
For $q$ a nondegenerate quadratic form over
$F$, recall that there is a short exact sequence of algebraic groups  
\begin{neqn} \label{Spin.seq}
\begin{CD}
1 @>>> C @>>> \Spin(q) @>>> SO(q) @>>> 1
\end{CD}
\end{neqn}
with $C \cong \mmu2$.

\begin{lem} \label{spinlem} \label{structlem}
For $q$ a $d$-dimensional nondegenerate quadratic form with
anisotropic part of  
dimension $d_{\mathrm{an}}$ such that $d \ge 5$ and 
$d + d_{\mathrm{an}} < 16$, the kernel of the Rost invariant of
$\Spin(q)$ is precisely the image of $H^1(F, C)$ in $H^1(F, \Spin(q))$.
\end{lem}

\begin{pf} 
The set $H^1(F, SO(q))$ classifies quadratic forms of the same
dimension and discriminant as $q$ \cite[29.29]{KMRT}.  For 
$\alpha \in H^1(F, \Spin(q))$ 
we set $q_\alpha$ to be the quadratic form corresponding to the image of
$\alpha$ in $H^1(F, SO(q))$.  Then $q_\alpha - q$ is not only
even-dimensional with trivial discriminant (i.e., $q_\alpha - q \in \IF2$),
but since $q_\alpha$ comes from $H^1(F, \Spin(q))$, it has the same
Clifford invariant as $q$ \cite[31.11]{KMRT} and so $q_\alpha - q \in \IF3$
by Merkurjev's Theorem. 
As described in \cite[p.~437]{KMRT}, the Rost invariant of $\alpha$ is
the Arason invariant $e_3(q_\alpha - q) \in H^3(F, \Z/2)$.  (Since $\Z/2 =
\mmud{2}$, we can consider $\Z/2$ to be a subgroup of $\QZt$ and hence
$H^3(F, \Z/2)$ is a subgroup of $H^3(F, \QZt)$.)	

Suppose first that $\alpha$ is in the image of $H^1(F, C)$.  Sequence
\eqref{Spin.seq} induces an exact sequence
\begin{neqn} \label{Spin.coseq}
\begin{CD}
SO(q)(F) @>>> H^1(F, C) @>>> H^1(F, \Spin(q)) @>>> H^1(F, SO(q)),
\end{CD}
\end{neqn}
and since the Rost invariant $R_{\Spin(q)}$ ``factors through''
$H^1(F, SO(q))$, certainly $R_{\Spin(q)}(\alpha)$ is trivial.

Conversely, suppose that $\alpha$ is in the kernel of the Rost invariant.
Then $e_3(q_\alpha - q)$ is trivial, but as mentioned in the
introduction the kernel of $e_3$ is precisely $\IF4$.
Since $\dim q_\alpha = \dim q = d$,
the hypotheses on $q$ ensure that the dimension of the anisotropic
part of $q_\alpha - q$ is strictly less than 16.  Since $q_\alpha - q
\in \IF4$, 
it is hyperbolic by the Arason-Pfister Hauptsatz \cite[X.3.1]{Lam}.
Thus $q_\alpha \simeq q$ and $\alpha$ is in the kernel of the map $H^1(F,
\Spin(q)) \ra H^1(F, SO(q))$, which is just the image of $H^1(F, C)$.
\end{pf}

The first map in \eqref{Spin.coseq} is the spinor norm, which
immediately produces the following lemma.

\begin{cor}\label{sncor}
Suppose that $q$ is as in Lemma \ref{spinlem}.  Then the kernel of
the Rost invariant is isomorphic to $\Fx / SN(q) \Fxsq$, where $SN(q)$ is
the image of the spinor norm map $SO(q)(F) \ra \Fsq$.\hfill $\qed$
\end{cor}

\begin{borel*} \label{Bcor}
Quasi-split simply connected groups of type $B_n$ are actually split,
so of the form 
$\Spin(q)$ for $q = n\H \perp \qform{1}$.  In terms of the lemma, $d =
2n + 1$ and $d_\mathrm{an} = 1$.  So $q$ satisfies the hypotheses for $2 \le n
\le 6$.  Since $q$ is isotropic, it has surjective spinor norm, so the
Rost invariant for a split group of type $B_n$ has trivial kernel for
$2 \le n \le 6$.
\end{borel*}

\begin{eg}[$B_7$] \label{B7.eg}
As just mentioned, 
the split simply connected group of type $B_7$ is isomorphic to
$\Spin(q)$ for $q = 7\hyp \perp \qform{1}$.  We will show that the
Rost invariant $R_{\Spin(q)}$ can have nontrivial kernel.  Sequence
(\ref{Spin.seq})
induces an exact sequence
\begin{neqn} \label{B7.seq}
\begin{CD}
H^1(F, \Spin(q)) @>>> H^1(F, SO(q)) @>\partial>> H^2(F, \mmu2)
\end{CD}
\end{neqn}
where the set $H^1(F, SO(q))$ classifies nondegenerate quadratic forms with
the same dimension (15) and discriminant $(1\cdot \Fxsq)$ as $q$.

Fix a base field $F$ and a nonhyperbolic 4-fold Pfister form $\varphi$ over $F$
(e.g. $F = \R$, $\varphi = \pform{-1, -1, -1, -1}$).
Set $q_\alpha = -\varphi'$ for $\varphi'$ such that $\varphi = \qform{1}
\perp \varphi'$.  Then $\disc q_\alpha = (-1)^{\binom{15}{2}}
\det(-\varphi') = 1 \cdot \Fxsq$, so there is a unique element of
$H^1(F, SO(q))$ corresponding to $q_\alpha$.  The image of $q_\alpha$
under the connecting homomorphism $\partial$ is $[C_0(q_\alpha - q)]$,
which by \cite[V.2.10]{Lam} is the same as $[C(q_\alpha - q)]$ which is
trivial since $q_\alpha - q = -\varphi \in \IF3$.  Thus $q_\alpha$ is the
image of some $\alpha$ in $H^1(F, \Spin(q))$.  But then
$R_{\Spin(q)}(\alpha) = e_3(q_\alpha - q) = e_3(-\varphi)$, which is
trivial since $\varphi \in \IF4$.
\end{eg}

\begin{borel*} \label{Dcor}
An analysis for groups of type $D_n$ similar to the one in \ref{Bcor} for
$B_n$ shows that the 
Rost invariant for a simply connected group is
trivial for groups of type $\oDn$ with $3 \le n \le 7$ and for groups
of type $\tDn$ with $3 \le n \le 6$.  As in the $B$ case, we show that
one of these bounds is sharp.
\end{borel*}

\begin{eg}[$\oDe$] \label{1D8.eg}
The situation here is quite similar to the one in Example \ref{B7.eg},
except that $q = 8\hyp$.  We use the same base field $F$ and nonsplit
4-fold Pfister form $\varphi$ from before.  Then there is a unique
element of $H^1(F, SO(q))$ corresponding to $\varphi$ and since
$\varphi = \varphi - q \in \IF4$, the same reasoning shows that there
is a nontrivial class in $H^1(F, \Spin(q))$ which is the inverse image
of $\varphi$ and which has trivial Rost invariant.
\end{eg}

Lemma \ref{spinlem} easily deals with quasi-split groups of type
$\tAn$ of low rank.

\begin{cor} \label{Acor}
If $G$ is a quasi-split simply connected group of type
$\tAn$ with $n \le 5$, the kernel of the Rost invariant $R_G$ is trivial.
\end{cor}

\begin{pf}
Set $K$ to be the quadratic field extension of $F$ which splits
$G$ and take $(V,h^d)$ to be a ``maximally split'' $(n+1)$-dimensional
hermitian form over 
$K$.  (See below for a more explicit description.)  Then $G$ is none
other than $SU(V,h^d)$, the algebraic group 
with $F$-points
\[
SU(V, h^d)(F) = \{ g \in GL(V)(K) \mid 
\text{$h(gv, gv') = h(v,v')$ for all $v, v' \in V$ and $\det g = 1$}
\}.
\]
The {\em trace form} of $h^d$ is
defined to be the quadratic form $q^d$ on $V$ considered as a
$2(n+1)$-dimensional vector space over $F$ given by $q^d(v) =
h^d(v,v)$.  Then
\[
h^d = 
\begin{cases}
m\hyp&\text{if $n + 1 = 2m$,}\\
m\hyp \perp \qform{1}&\text{if $n + 1 = 2m + 1$}
\end{cases}
\qquad\text{and}\qquad
q^d = 
\begin{cases}
2m\hyp&\text{if $n + 1 = 2m$,}\\
2m\hyp \perp \pform{d}&\text{if $n + 1 = 2m + 1$,}
\end{cases}
\]
where $K = F(\sqrt{d})$ if $n = 2m$ for some integer $m$, and the
$\hyp$ occurring in the description of $h^d$ is the usual unitary
hyperbolic plane as described in \cite[7.7.3]{Sch}.

The set $H^1(F, G)$ classifies nonsingular hermitian forms $h$ on $V$
which have the same dimension and discriminant as $h^d$
\cite[p.~403]{KMRT}. 
Now $G \injects SO(V, q^d)$ and the corresponding map $H^1(F, G) \ra
H^1(F, SO(V, q^d))$ sends $h$ to its trace form $q$, and this map is an
injection by \cite[10.1.1(ii)]{Sch}.  Moreover, the Rost invariant
$R_G(h)$ is just $e_3(q - q^d)$ by \cite[31.44]{KMRT}.  Since $\dim q^d
= 2n + 2 < 14$ and the anisotropic part of $q^d$ has dimension 0 (if
$n + 1$ is even) and 2 (if $n + 1$ is odd), as in the proof of Lemma
\ref{structlem}, if $R_G(h)$ is trivial, $q \simeq q^d$ and so $h
\simeq h^d$.
\end{pf}

\begin{eg}[$\tAs$] \label{2A6.eg}
Take $F = \R$, $K = \mathbb{C}$, and consider $G = SU(V, h^d)$ for $h^d$ the
hermitian form $3\hyp \perp \qform{1}$ over $K$, so that $G$ is simply
connected quasi-split of type $\tAs$.  Then the hermitian form $h =
\qform{-1, -1, -1, -1, -1, -1, -1}$ has trace form $q = -7\pform{-1}$
which is not hyperbolic, so $h$ corresponds to a (unique) nontrivial
class in $H^1(F, G)$.  However,
\[
q - q^d = -7\pform{-1} - \pform{-1} = -\pform{-1, -1, -1, -1} \in
\IF4,
\]
so $R_G(h)$ is trivial.
\end{eg}

\section{Folded root systems} \label{foldsec}

\begin{borel}{The Rost multiplier} \label{Rm.def}
In general, for an arbitrary algebraic group $G$ we define the set $G_\ast$ of
loops in $G$ to be the homomorphisms $\mathbb{G}_m \ra G$.  Then as in
\cite[p.~432]{KMRT}, we
set $Q(G)$ to be the abelian group of all integer-valued functions on
$G_\ast$ such that
\begin{enumerate}
\item for $^gf$ the loop given by $^gf(x) = gf(x)g^{-1}$, $q(^gf) = q(f)$ for all $g \in G$ and $f \in G_\ast$; and
\item for any two loops $f$ and $h$ with commuting images, the
function $\Z \times \Z \longto \Z$ given by $(k, m) \mapsto q(f^k
h^m)$ is a quadratic form.
\end{enumerate}

When $G$ is an almost simple group, $Q(G)$ is cyclic with a canonical
generator which is positive definite \cite[31.27]{KMRT}, hence is
identified with $\Z$.  
Now suppose that we have two almost simple simply connected groups $H
\injects G$.  This map clearly
induces a map $H_\ast \ra G_\ast$, so we in turn have a map $\Z = Q(G) \ra
Q(H) = \Z$.  Because the canonical generators are positive definite,
this map must be multiplication by a positive integer $n$, which 
we define to be the {\em Rost multiplier} of the embedding.

The naturality of the Rost invariant implies that 
we have a commutative diagram
\[
\begin{CD}
H^1(F, H) @>>{R_H}> H^3(F, \QZt) \\
@VVV @V{n\cdot}VV \\
H^1(F, G) @>{R_G}>> H^3(F, \QZt),
\end{CD}
\]
where $n$ is the Rost multiplier of the embedding \cite[31.34]{KMRT}.
This is the motivation for our study of this invariant.
\end{borel}

\begin{borel*} \label{Rm.comp}
Luckily, it can be quite easy to compute such a ``Rost
multiplier''.  Suppose that $G$ and $H$ are split and contain split
maximal tori $S$ and $T$ respectively such that the image of $T$ lies
in $S$.
Since $G$ and $H$ are simply connected, the character groups
$X(T)$ and $X(S)$ are identified with the weight lattices, but the
character groups are dual to the loop groups $S_\ast$ and $T_\ast$
\cite[8.6]{Borel} and
the weight lattices are dual to the lattices generated by the coroots,
which we denote by $\Lambda_{c,G}$ and $\Lambda_{c,H}$, respectively.
(By a {\em coroot}, we mean the roots of the dual root system, which
are denoted by $\ach$ in \cite[VI.1]{Bou:g4} for $\alpha$ a root.)
Putting these dualities together, we obtain identifications $S_\ast =
\Lambda_{c,G}$ and $T_\ast = \Lambda_{c,H}$, so 
the embedding $T \ra S$ induces a map $\Lambda_{c,H} \ra
\Lambda_{c,G}$.  Now the dual root systems (whose roots are the
coroots) are indeed root systems \cite[VI.1.1, Prop.~2]{Bou:g4} and so
they each 
have a unique minimal Weyl-group invariant positive-definite
integer-valued quadratic form \cite[VI.1.2, Prop.~7]{Bou:g4}, say $q$
and $r$ (for  
the forms for $G$ and $H$ respectively).  Hence $q$ induces such a
form on $\Lambda_{c,H}$, which must be of the form $nr$ for some
natural number $n$.  This $n$ is the Rost multiplier of the inclusion.

Criterion (2) in the definition of $Q(G)$ implies that its canonical
generator is identified with the positive-definite Weyl-group
invariant quadratic form on the dual root  
system which takes the value 1 on short coroots.  (Short roots correspond 
to long roots, where we adopt the convention that short = long in the
event that all roots have the same length.  In that case, the
quadratic form is very easy to identify, in that its Gram matrix is
simply the Cartan matrix of the root system with all entries divided
by 2.)  So one can simply compute the image of a short coroot from
$H$ in the dual root system for $G$ to find the Rost multiplier of the
embedding.
\end{borel*}

\begin{eg}[$SL_n \ra SL_{2n}$]
The block diagonal embedding $SL_n \injects SL_{2n}$ via $x \mapsto
\stbtmat{x}{}{}{x}$ has Rost multiplier 2.  The embedding given by $x
\mapsto \stbtmat{x}{}{}{1}$ has Rost multiplier 1.
\end{eg}

\begin{eg}[Folding] \label{folding} \label{E6.rm} \label{e6rm}
The split simply connected group of type $E_6$ can be realized as the
group $\Inv(J)$ of invertible linear maps of the split Albert algebra
$J$ which preserve the cubic norm form.  The algebra $J$ has a
nondegenerate symmetric bilinear trace form $T$ given by setting
$T(x,y)$ to be the trace of the product $x \cdot y$ \cite[p.~240,
Thm.~5]{Jac:J}, and for $\varphi 
\in \Inv(J)(F)$ we define $\varphi^\dag \in GL(J)(F)$ to be the unique map
satisfying $T(\varphi(j), \varphi^\dag(j')) = T(j,j')$ for all $j, j'
\in J$.  This defines an outer automorphism of $E_6 = \Inv(J)$
\cite[p.~76, Prop.~3]{Jac:J3} and the subgroup of elements fixed by
this automorphism 
is the split group $F_4$ of $F$-algebra automorphisms of $J$.

We would like to compute the Rost multiplier of the embedding $F_4
\subset E_6$.  We fix an $F$-split maximal torus $S$ in $G := E_6$
which is preserved by the automorphism (such as the one denoted by
``$S_6$'' in \cite[pf.~of 7.2]{G:struct}) and 
fix a set of simple roots $\D$ of $G$ with respect to $S$.  We would
like our outer automorphism to leave $\D$ invariant, although it
probably does not do so.
However, two things are apparent from the definition of the Rost
multiplier: it is not changed if we extend scalars nor if modify 
the automorphism $\varphi \mapsto \varphi^\dag$ by an inner
automorphism of $E_6$.  So we may assume that the base field is
separably closed 
and so that the $F$-points of the Weyl group of $G$ with respect to
$S$ (i.e., the $F$-points of $N_G(S) / S$) is the full Weyl group of
the root system of $G$ with respect to $S$.  Then we may modify our
outer automorphism by an element of the Weyl group so that $F_4$ is
described as the subgroup of $E_6$ fixed by the automorphism $f$
induced by the automorphism of $\D$ which is given by the unique
nontrivial automorphism of the Dynkin diagram.  That is, we set $H :=
F_4 = G^f$ (= the subgroup of $G$ of elements fixed by $f$), and $T :=
\cc{(S^f)}$ (= the identity component of $T \cap G^f$) is a maximal
torus in $H$.  Then the restrictions of elements of $\D$ to $T$ give a
root system of $H = F_4$ with respect to $T$ \cite[p.~108]{Schatt} and
the fibers of this restriction map are the orbits of $f$ in $\D$
\cite[3.5]{Schatt}.

Now $\Lambda_{c,G}$ is a free $\Z$-module with basis $\Dch = \{ \dch
\mid \delta \in \D \}$ which is permuted by $f$ and $\Lambda_{c,H}$ is
the fixed sublattice.  So $\Lambda_{c,H}$ has a basis consisting of
one element for each orbit of $f$ in $\Dch$,
and this element is given by the sum of the elements in the orbit in $\Dch$.
Since no two roots in $\D$ lying in the same orbit under $f$ are
connected in the Dynkin diagram, the description of the roots of $H$
with respect to $T$ 
above shows that the spanning set in $\Lambda_{c,H}$ just described is
in fact a set of simple coroots of $H$ with respect to $T$, which is
dual to the set of simple roots given by restrictions of elements of
$\D$.

So we would like to compute the value of $q$ on a particular simple coroot
$\ach$ of $H$ with respect to $T$ 
under the composition $\Lambda_{c,H} \ra \Lambda_{c,G}
\xrightarrow{q} \Z$.
Since each orbit of $f$ in $\Dch$ consists
of totally disconnected sets of vertices of the Dynkin diagram and all
coroots $\dch$ in $\Dch$ have $q(\dch) = 1$ as a consequence of
hypothesis (2) in the definition of $Q(G)$, the value of $q$ on $\ach$
is precisely the size of the orbit in $\Dch$ which gave rise to
$\ach$.  Since there is some element of $\Dch$ is fixed by $f$, $q$ takes the
value 1 on some $\ach$, so by the discussion in \ref{Rm.comp} the
Rost multiplier of the embedding $H \subset G$ is one.
\end{eg}

\begin{rmk}
Presumably this same argument also works in the other instances where
one obtains a root system by ``folding up'' another root system all of
whose roots have the same length, i.e., $C_{\ell + 1} \subset A_{2\ell
+ 1}$, $B_{n-1} \subset D_n$, and $G_2 \subset D_4$.  The other root
system consisting of roots of the same length, $A_{2\ell}$, folds up
to give the smaller root system $BC_\ell$, see \cite[Table I]{Heck}.
\end{rmk}
\section{Small representations} \label{smallsec}

We say a representation $V$ of an algebraic group $G$ is {\em small}
if $G$ has an open orbit in $\pr(V)$.  We are interested in small
representations in the case where $G$ is simple, which have all been
classified as a consequence of the (more general) classification of
prehomogeneous vector spaces, see \cite{Kim:surv} for a survey.  
These small representations also provide ``standard relative
sections'' in the language of \cite[1.7]{Pop:sec}, and in that sense
were classified in \cite[Table 1]{Elash}.
Our motivation for studying these representations comes from the
following easy lemma, which was pointed out to me by Rost.

\begin{lem} \label{cosurj}
Suppose that $G$ is an
algebraic group over a field $F$ such that $G$ has a small
representation $V$, and 
that $F$ is infinite or $G$ is connected.  Let $H$ be 
the subgroup of $G$ consisting of elements which stabilize 
some $F$-point in the open orbit in $\pr(V)$.  Then the natural map
\[
H^1(F, H) \ra H^1(F, G)
\]
is surjective.
\end{lem}

\begin{pf}
If the base field $F$ is finite, then by hypothesis $G$ is connected, and by Lang's
Theorem $H^1(F, G)$ is trivial so the lemma holds.  So we may assume
that $F$ is infinite. 

Fix a 1-cocycle $z \in Z^1(F, G)$.  For $U$ the open orbit in
$\pr(V)$, $z$ defines a new open subset $U_z$ in $\pr(V)$ which is the
same as $U$ over the separable closure $F_s$ of $F$ but has a different 
Galois action: For $u \in
U_z(F_s)$ and $\s \in \Gal(F_s/F)$, $\s$ acts by
\[
\s \ast u = z_\s \s u
\]
where juxtaposition denotes the usual action.

Since the representation gives a map $G \ra GL(V)$, $\pr(V)_z \cong
\pr(V)$.  Thus since $F$ is infinite, $\pr(V)_z(F)$ is dense in
$\pr(V)_z(F_s)$.  Since $U_z(F_s)$ is open in $\pr(V)_z(F_s)$, the two
sets $U_z(F_s)$ and $\pr(V)_z(F)$ must meet nontrivially, i.e., 
$U_z$ has some $F$-point which we will denote by $x_z$.

Now let $x \in U(F)$ be the point with stabilizer subgroup $H$ and fix
some $g \in G(F_s)$ such that $gx = x_z$.  Then for all $\s \in
\Gal(F_s/ F)$, the element $g^{-1} z_\s \s g \s^{-1}$ fixes $x$ and so
lies in $H(F_s)$.  Thus $z$ is cohomologous to something in the image
of $Z^1(F, H)$.
\end{pf}

\begin{eg}[$O_{n-1} \subset O_n$]
Write $O_n$ for the orthogonal group of the dot product on
$F^n$.  Then the subgroup of $O_n$ which stabilizes $[v] \in
\pr(F^n)$ where $v$ has nonzero length is just $O_{n-1} \times \mmu2$,
where $O_{n-1}$ is the orthogonal group for the $(n-1)$-dimensional
space of vectors in $F^n$ which are orthogonal to $v$.  Iterating this
process recovers the fact that all nondegenerate quadratic forms are
diagonalizable, a.k.a.~the Spectral Theorem.
\end{eg}

\begin{eg}[$\Spin_n$ \cite{Igusa}, \cite{GV}, \cite{Pop:14}]
For $\Spin_n$ the spin group for an $n$-dimensional maximally split
quadratic form, the spin representation (if $n$ is odd) or the
half-spin representation (if $n$ is even) is small
for $n \le 12$ and $n = 14$.  In the $n = 14$ case, the stabilizer
subgroup is isomorphic to $(G_2 \times G_2) \semid \mmu8$, and this
leads to structural statements about 14-dimensional forms in $\IF3$,
see \cite{Rost:14}.
\end{eg}

\begin{eg}[$F_4 \times \mmu3 \subset E_6$] \label{e6sm}
We write $E_6$ for the split group of type $E_6$ which can be realized
as $\Inv(J)$ as described in \ref{E6.rm}. 
By \cite[p.~71, Thm.~7]{Jac:J3}, 
$E_6$ acts transitively 
on the subset of $J$ consisting of elements of norm 1, so certainly
this is a small representation.

We take $H$ to be the subgroup of $E_6$ which fixes the identity
element $1_J$ of $J$ projectively.  
Then since the norm form is cubic,
$\mmu3$ is contained in $H$ and is central (since it consists of scalar
endomorphisms), and any element of $H$ differs by an
element of $\mmu3$ from something which fixes $1_J$ absolutely.  This
subgroup of elements fixing $1_J$ is well-known --- it is none other
than the automorphism group  
$F_4$ of $J$ \cite[p.~186, Thm.~4]{Jac:J1}, which is split of type
$F_4$.  So $H \cong F_4 \times \mmu3$, and the resulting 
surjective map $H^1(F, F_4 \times \mmu3) \ra H^1(F, E_6)$ is the
statement that $H^1(F, E_6)$ classifies cubic forms of the form  
$\la N$ for $N$ the norm form on some Albert $F$-algebra and $\la \in
\Fx$,  see \cite{Sp:cubic}.
This can also be interpreted in terms of structurable algebras, see
\cite[2.8(1)]{G:struct}.
\end{eg} 

\begin{eg}[$E_6 \semid \mmu4 \subset E_7$]
\label{e7sm} \label{E7.sm}
Write $E_7$ for the split simply connected group of type $E_7$ over $F$.
It is the group of vector space automorphisms of $V = \stbtmat{F}{J}{J}{F}$
which preserve a quartic form $q$ as given in \cite[p.~87]{Brown:E7}.
Then $E_7$ acts transitively on the open subset of $\pr(V)$ consisting 
of points $[v]$ such that $q(v) \ne 0$ by \cite[7.7]{Ferr:strict}.

We set $H$ to be the subgroup of $E_7$ which stabilizes the vector 
$v = \stbtmat{1}{0}{0}{1}$ projectively.  This vector has $q(v) \ne
0$, and so by \cite[3.7]{Ferr:strict} there are two uniquely
determined (up to scalar multiples) ``strictly regular'' elements
$e_1$ and $e_2$ such that $v$ lies in their span.  These are $e_1 =
\stbtmat{1}{0}{0}{0}$ and $e_2 = \stbtmat{0}{0}{0}{1}$.  Since $E_7$
preserves the property of being strictly regular, every element of $H$
must projectively stabilize $e_1$ and $e_2$ as well, and perhaps interchange
them.

Now, the map $\omega$ defined by
\[
\omega \basmat = \tbtmat{i \beta}{i j'}{i j}{i \alpha}
\]
lies in $H$, where $i$ is some fixed square root of $-1$ in the
separable closure of $F$.  We would like to describe an arbitrary $h
\in H$, which after modification by $\omega$ we may assume projectively
stabilizes each of $e_1$ and $e_2$.  Then by \cite[p.~96,
Lem.~12]{Brown:E7}, $h$ must be of the form
\[
h \basmat = \tbtmat{\mu^{-1}
\alpha}{\varphi(j)}{\varphi^\dagger(j')}{\mu \beta}
\]
where $\varphi$ is a similarity of the norm form on $J$ with
multiplier $\mu$ and $\varphi^\dagger$ is as defined in \ref{e6rm}.
Since $h$ also stabilizes $v$, we must have that
$\mu = \pm 1$.  In particular, after modifying $h$ by $\omega^2 = -1$,
we may assume that $h$ has the form
\[
h \basmat = \tbtmat{\alpha}{\varphi(j)}{\varphi^\dagger(j')}{\beta}
\]
where $\varphi$ preserves the cubic norm on $J$ and so lies in $E_6$.
So we have shown that $H \cong E_6 \semid \mmu4$.

The surjection on Galois cohomology coming from this example will be
more useful if we can replace $E_6 \semid \mmu4$ with an almost simple
group.  For $K$ a quadratic \'etale $F$-algebra, we write $E_6^K$ for
the simply connected quasi-split group of type $E_6$ over $F$ which is
split by an extension $L$ of $F$ if and only if $L \otimes_F K \cong L
\times L$.

\begin{prop} {\rm (Cf.~\cite[4.14]{G:struct})} \label{e7smprop}
Suppose that $F$ has characteristic $\ne 2, 3$.
For each $\alpha \in H^1(F, E_7)$ there is some quadratic \'etale
$F$-algebra $K$ such that $E_6^K$ embeds in $E_7$ with Rost multiplier
$1$ and $\alpha$ is in
the image of the induced map $H^1(F, E_6^K) \ra H^1(F, E_7)$.
\end{prop}

\begin{pf}
Fix some $a \in Z^1(F, E_6 \semid \mmu4)$ representing $\alpha$.
The natural projection $E_6 \semid \mmu4 \ra \mmu4$ has an obvious
section given by sending $i \mapsto \omega$, and we set $b$ to be the
image of $a$ given by the map induced by the composition $E_6 \semid
\mmu4 \ra \mmu4 \ra E_6 \semid \mmu4$.  We twist $E_6 \semid \mmu4$ by
$b$ to obtain a new group $(E_6 \semid \mmu4)_b$, with a twisted
Galois action $*$ so that 
\[
\s * g = b_\s ({^\s g}) b_\s^{-1},
\]
where ${^\s g}$ denotes the usual action.  Then we have an 
isomorphism
\[
\begin{CD}
H^1(F, (E_6 \semid \mmu4)_b) @>{\sim}>{\tau_b}>
H^1(F, E_6 \semid \mmu4)
\end{CD}
\]
where $\tau_b^{-1}(\alpha)$ is the class of a 1-cocycle given by $\s
\mapsto a_\s b_\s^{-1}$ with values in the 
identity component of the twisted group $(E_6 \semid
\mmu4)_{b}$.  This identity component is just $E_6$ twisted by
$b$, and we would like to show that it is isomorphic to $E_6^K$ for some
quadratic \'etale $F$-algebra $K$.  We observe that if $\s$ in the
absolute Galois group of $F$ has $b_\s = \pm 1$, then $\s$ acts in the
usual manner upon the twisted $E_6$.  On the other hand, if $b_\s =
\pm \omega$, then the twisted action is given by
\[
(\s * h)\basmat = (\pm \omega) \s h \s^{-1} (\pm \omega)^{-1} \basmat
= \tbtmat{\alpha}{\s \varphi^\dagger
\s^{-1}(j')}{\s\varphi\s^{-1}(j)}{\beta}.
\]
Since this is precisely the description of the Galois action on
$E_6^K$ given in \cite[2.4]{G:struct} for $K$ determined by the image
of $b$ under the composition $H^1(F, E_6 \semid \mmu4) \ra H^1(F,
\mmu4) \ra H^1(F, \mmu2) = \Fsq$, we have $(E_6)_b \cong E_6^K$. 
To see that $E_6^K$ embeds in $E_7$, we observe that the 1-cocycle $b$ is 
trivial in $H^1(F, E_7)$ by \cite[4.10, 5.10]{G:struct}, so we have a map
\[
\begin{CD}
E_6^K \subset (E_6 \semid \mmu4)_b \injects (E_7)_b @>{\sim}>f> E_7
\end{CD}
\]
where (by a simple computation having nothing to do with $E_7$)
$H^1(f) = \tau_b$.  This proves the proposition aside  
from the claim about the Rost multiplier.

But that claim is easy in the split case (where $K = F \times F$),
since the embedding of $E_6$ in $E_7$ comes from the obvious embedding
of root systems.
Then since the Rost multiplier is invariant under scalar extensions, the
embeddings of quasi-split groups of type $E_6^K$ in $E_7$ given above
all have Rost multiplier 1 as well.
\end{pf}
\end{eg}

\section{$\agrp{1}{D_4} \subset \dEvi$} \label{Gsec}

For the remainder of the paper we will study the quasi-split group
$E^K_6$ of type $\dEvi$ defined in \ref{e7sm}.  In this section we
introduce a particular subgroup $G$ of $E^K_6$ which is reductive of
semisimple type $\oD$.  Defining $G$ will necessitate digging more
deeply in to the structure of Albert and Cayley algebras.

\begin{defn}
Fix $\C$ to be the split Cayley algebra endowed with hyperbolic norm
form $\n$ and canonical involution denoted by $\bar{\ }$.  We define
some algebraic groups related to $\C$.  (For more information about
Cayley algebras, see \cite[\S 33.C]{KMRT} or \cite[Ch.~III, \S
4]{Schfr}.)  First, if $t \in GL(\C)$
satisfies $\n(t(c)) = m \n(c)$ for some $m \in \Fx$ and all $c \in
\C$, we say that $m$ is a {\em similarity} of $\n$ with multiplier
$\mu(t) := m$.  (Note that if $\s_\n$ is the involution on $\EndF(\C)$
which is adjoint for $\n$ so that $\n(tc, c') = \n(c, \s_\n(t) c')$
for all $c, c' \in \C$, then $\mu(t) = \s_\n(t) t$.)  
Then we set $\GOP\Cn$ to be the algebraic group with
$F$-points
\[
\GOP\Cn(F) := \left\{ t \in GL(\C) \left| \parbox{3in}{$t$ is a similarity
of $\n$ with multiplier $\mu(t)$ such that $\det(t) = \mu(t)^4$}
\right. \right\}.
\]

We can also
define a new, seemingly uglier multiplication $\star$ on $\C$ by
setting $x \star y := \bar{x} \bar{y}$ as in \cite[\S 34.A]{KMRT}.  
Then
a {\em related triple} is a triple $\trion{t}$ in $\trio{\GOP\Cn}$
such that 
\[
\mu(t_i)^{-1} t_i(x \star y) = t_{i+2}(x) \star t_{i+1}(y)
\]
for all $x, y \in \C$ and $i = 0, 1, 2$ with subscripts taken modulo
3.
We write $\Rel\Cn$ for the algebraic subgroup of $\trio{\GOP\Cn}$
consisting of related triples and  
$\Spin(\n)$ for the subgroup of $\Rel\Cn$ consisting of triples with
multiplier one 
(i.e., those triples such that $\mu(t_i) = 1$ for all $i$).
\end{defn}
\begin{borel*}
The vector space underlying the split Albert $F$-algebra $J$ is the
subspace of $M_3(\C)$ 
consisting of elements fixed by the conjugate transpose $\ast$ which
applies $\bar{\ }$ to each entry and takes the transpose.  It is the
algebra denoted by $\mathfrak{\h}(\C_3)$ in the notation of \cite[\S
I.5]{Jac:J} and has multiplication $a \cdot b := (ab + ba) /2$, where
juxtaposition denotes the usual multiplication on $M_3(\C)$.
When
writing down explicit elements of $J$, we will use a ``$\cdot$'' to
indicate entries whose values are forced by this symmetry condition.
Then $\Rel\Cn$  embeds in the group $\Inv(J)$ of norm isometries of
$J$ via the map $\ut \mapsto 
 g_\ut$ given by  
\begin{neqn} \label{gdef}
g_\ut \basjord = \jordmat{\mu(t_0)^{-1} \e_0}{\mu(t_1)^{-1}
\e_1}{\mu(t_2)^{-1} \e_2}{t_0(c_0)}{t_1 (c_1)}{t_2 (c_2)}.
\end{neqn}
\end{borel*}

\begin{borel}{Definition of $G$}
Since $\Rel\Cn$ embeds in $\Inv(J)$ over $F$, it embeds in $E^K_6$
over $K$.  However, we can identify $E^K_6$ with
$\Inv(J)$ with a different $\iota$-action where $^\iota f := \iota
f^\dag \iota$, where $\iota$ is the nontrivial
$F$-automorphism of $K$ and juxtaposition denotes the usual action,
and we fix this identification for the rest of the paper.
Then the map $\Rel\Cn \ra E^K_6$ is not
defined over $F$, because  
for $\ut = \trion{t} \in \Rel\Cn(F_s)$ and $g_\ut \in E^K_6$, we have 
$^\iota g_\ut = g_{\iota \s_\n(\ut)^{-1} \iota}$ which is typically not
the same as $g_{\iota \ut \iota}$ where $\s_\n(\ut)$ means to apply
$\s_\n$ to each component of $\ut$.
So we define $G$ to be the algebraic group over $F$
which is the same as $\Rel\Cn$ over $K$ but with a different
$\iota$-action: for $\ut \in G(F_s)$, we set $^\iota \ut := \iota
\s_\n(\ut)^{-1} \iota$.  Then $G$ injects into $E^K_6$ over $F$ via
the map $g$ from \eqref{gdef}.
\end{borel}

This group $G$ is
reductive with absolute rank 6 and semisimple part $\Spin(\n)$ of type
$\oD$.  Its importance is 
given by the following lemma, excavated from a paper by Ferrar:

\begin{ferrlem} \label{ferrlem} {\rm \cite[p.~65, Lem.~3]{Ferr:E6}} 
The natural map $H^1(K/F, G) \ra H^1(K/F, E_6^K)$ is surjective.
\end{ferrlem}

\begin{comm}
Ferrar proved this by explicit computations in the Jordan algebra.
However, this can also be seen with more algebraic group-theoretic
methods, as was pointed out to me by Gille: We must assume that our
base field has characteristic zero, which as we observed in the
introduction doesn't impair our main results in any way.  Let $T$ be a
maximal torus in $G$ defined over $F$.  Then \cite[p.~329,
Lem.~6.17]{PlatRap} combined with the fact that all maximal tori in
$E_6^K$ are conjugate over $F$ says that there is a Borel subgroup $B$
of $E_6^K$ defined over $K$ such that $T = B \cap {^\iota\hskip-.1em
B}$.  Consequently, the natural map $H^1(K/F, T) \ra H^1(K/F, E^K_6)$
is a surjection by \cite[p.~369, Lem.~6.28]{PlatRap}. 
\end{comm}

Now imagine how the argument for proving the main theorem in the $\dEvi$
case must proceed: We apply some simple argumentation and Ferrar's
Lemma to show that any class in $H^1(F, E^K_6)$ with trivial Rost
invariant must come from $H^1(K/F, G)$.  Then we apply some facts
about Rost invariants on this smaller group to obtain the theorem.
However, $G$ is reductive, so we want to put our class
with trivial Rost invariant into a semisimple subgroup if we hope to
apply our results from Section \ref{BDsec}.  So we need to do
something about the center of $G$.

\begin{borel}{The center $P$ of $G$}
For the moment, set $N_1$ to be the algebraic group with $F$-points
the elements of $\Kx$ with norm 1 in $F$.  This group is the same as
$\Gm$ over $K$, but has a different $\iota$-action given by $^\iota
\lambda = \iota(\lambda)^{-1}$.  It is often denoted by something like
$R^{(1)}_{K/F}(\mathbb{G}_{m,K})$. 

Now, since $\Gm$ is the center of $\GOP\Cn$, the center of $\Rel\Cn$
is the subgroup of $\trio{\Gm}$
consisting of triples whose product is one.  But we are concerned with
$G$, which has a different $\iota$-action; its center $P$ is then
isomorphic to the subgroup of $N_1^{\times 3}$ consisting of triples
whose product is 1.
\end{borel}

\begin{borel*}\label{Pdesc}
We have a map $G \ra P$ given by sending each $t_i$
to its multiplier $\mu(t_i) = \s_\n(t_i) t_i \in N_1$.  The semisimple part
$\Spin(\n)$ of $G$
is the kernel of this map, and we have a short exact sequence 
$1 \ra \Spin(\n) \ra G \ra P \ra 1$.
This sequence is even exact over $K$ (instead of just over a separable
closure of $F$) because the map $G \ra P$ is surjective over $K$ by
\cite[35.4]{KMRT}. 

So consider the map $H^1(K/F, G) \ra H^1(K/F, P)$.  We need to
describe this second group.  We have a short exact sequence over $K$
\[
\begin{CD}
1 @>>> P @>>> N_1^{\times 3} @>{\pi}>> N_1 @>>> 1,
\end{CD}
\]
where $\pi$ is the product map, which induces an exact sequence
\[
\begin{CD}
1 @>>> H^1(K/F, P) @>>> H^1(K/F, N_1^{\times 3}) @>{H^1(\pi)}>>
H^1(K/F, N_1).
\end{CD}
\]
(The second map is an injection because the product map $\pi$ is a
surjection on $F$-points.)  Now any
1-cocycle in $Z^1(K/F, N_1)$ is determined by its value at $\iota$,
and the condition that it is a 1-cocycle forces that this value lies
in $\Fx$.  Then the obvious check shows that two such are cohomologous
if and only if they differ by a norm from $\Kx$.  So $H^1(K/F, P)$ is
isomorphic to the subgroup of $\trio{(\Fx / N_{K/F}(\Kx))}$ with
product one.

So if we have a 1-cocycle $\gamma \in Z^1(K/F, G)$, it is determined
by its value at $\iota$, and the image of $\gamma$ in $H^1(K/F, P)$ is
the multiplier of its value there.
\end{borel*}

One natural question arising from this discussion is the
following: If you have a 1-cocycle in $H^1(K/F, 
E^K_6)$, then it comes from $H^1(K/F, G)$ by Ferrar's Lemma and so has
an image in $H^1(K/F, P)$.  Is that image an invariant of the original
class in $H^1(K/F, E^K_6)$?  The answer is a very definite no, as is
shown in the following lemma.  (We will give explicit situations
where the hypotheses are satisfied nontrivially in \ref{specialcor} and
\ref{closer.2E6}.)  Let $e_i \in J$ be the element whose
only nonzero entry is a 1 in the $(i+1, i+1)$-position.

\begin{movlem} \label{movlem}
Let $\eta$ be a 1-cocycle in $Z^1(K/F, G)$ whose image in $Z^1(K/F,
P)$ takes the value $\ua$ at $\iota$.  Suppose that there is some
$j \in e_0 \times J_K$ such that 
\[
j^\# = 0 \qquad \text{and} \qquad T(j, \eta_\iota \iota j) = r \in
\Fx.
\]
Then $\eta$ is cohomologous in $H^1(K/F, E^K_6)$ to a 1-cocycle coming from
$Z^1(K/F, G)$ whose image in $Z^1(K/F, P)$ takes the value 
$(r, a_0, (a_0 r)^{-1})$ at $\iota$.
\end{movlem} 

The hypotheses in the lemma make use of the Freudenthal cross product
$\times \!: J \times J \ra J$, which is a commutative bilinear map
defined by the relation
$6 N(j) = T(j, j \times j)$
for all $j \in J$.  The map $\# \!: J \ra J$ is defined by 
$2 j^\# := j \times j$.

\medskip
\begin{pf}
The proof is an adaptation of an argument in \cite[p.~277]{Ferr:class}.
We first observe that the
three elements $j$, $e_0$, and $e_0 \times j'$ for $j' := \eta_\iota
\iota j$ all have ``rank one'', i.e., are sent to zero by the map $x
\mapsto x^\#$.  Furthermore, since $e_0 \times J$ is preserved by
$\iota$ and $G(K)$ \cite[1.7]{G:struct}, $j' = e_0 \times y$
for some $y$ and
applying a linearization of \cite[(19)]{McC:FST} gives us
\[
e_0 \times (e_0 \times j') = T(e_0, e_0) (e_0 \times y) = j'.
\]
Thus this triple is ``in general position'', meaning that for $N$
trilinearized so that $N(x,x,x) = N(x)$, we have
\[
6 N(e_0, j, e_0 \times j') = T(e_0 \times j, e_0 \times j') = T(j,
e_0 \times (e_0 \times j')) = T(j, j') = r \ne 0.
\]
Thus by \cite[3.11]{SpV}, there is some $f \in
\Inv(J)(K)$ such that 
\[
f(j) = \rho_0 e_0, \quad f(e_0) = \rho_1 e_1, \ \text{and\ } \ 
f(e_0 \times j') = \rho_2 e_2
\]
for $\rho_i \in \Kx$.  Since $\n$ is hyperbolic, $G(K)$ is split over
$K$ and the map $G(K) \ra P(K)$ is surjective.  So there is some $g =
g_\ut \in \Inv(J)(K)$ such that $g(e_0) = \rho_0^{-1} e_0$, $g(e_1) =
\rho_1^{-1} e_1$, and $g(e_2) = \rho_0 \rho_1 e_2$.  By replacing $f$
with $gf$, we may assume that $\rho_0 = \rho_1 = 1$.
Moreover, $f$ preserves $N$, and so 
$r = 6 N(e_0, e_1, \rho_2 e_2) = \rho_2$.

So we set $\eta' \in Z^1(K/F, E^K_6)$ to be the cocycle cohomologous
to $\eta$ given by $\eta'_\iota = f^\dag \eta_\iota \iota f^{-1} \iota$.
Keeping in mind the facts that $e_i \times e_{i+1} = e_{i+2}$; $f(u
\times v) = f^\dag(u) \times f^\dag(v)$ for all $u, v \in J_K$; and
$j \times (e_0 \times j') = r e_0$ (as can be verified by
examining the explicit formula for $\times$ given in \cite[p.~358,
(4)]{Jac:J}, although the reader should be cautioned that our
definition of $\times$ --- which agrees with the one in \cite{KMRT}
and \cite{McC:FST} --- differs from Jacobson's by a factor of 2), one 
can now easily calculate that
$f^\dag(e_0) = e_1$ and $f^\dag(j') = r e_0$.
Then it follows that
\[
\eta'_\iota(e_0) = r e_0, \quad and \quad \eta'_\iota(e_1) = a_0^{-1}
e_1.
\]
Finally, since $\eta_\iota \iota (e_0 \times j') = a_0 e_0 \times j$,
we have $\eta'_\iota(e_2) = (a_0 /r)e_2$, and we are done.
\end{pf}

\section{$\agrp{2}{D_5} \subset \dEvi$}

For the purpose of making computations, we will need to make use of
another subgroup of $E^K_6$, which we define to be the subgroup
consisting of elements $h$ such that $h$ and $h^\dag$ both fix the
element
$e_0 \in J$.
Since the map $h \ra h^\dag$ is a group homomorphism on $\Inv(J)$, it
is clear that  
$H$ is indeed a subgroup of $E^K_6$ over $K$, and it is certainly preserved by
the $\iota$-action so it is even defined over $F$.  Our first task is to
describe it explicitly.

\begin{borel*} \label{udef}
Fix a
particular basis $u_1,
u_2, \ldots, u_8$ for the split Cayley algebra $\C$ as given in
\cite[p.~388]{G:iso}.  One important thing for us to know about this
basis is that when we bilinearize the norm form $\n$ so that $\n(x,x)
= \n(x)$, we have
\[
\n(u_i, u_j) = \begin{cases}
1 & \text{if $i + j = 9$} \\
0 & \text{otherwise},
\end{cases}
\]
so that the Gram matrix of $\n$ with respect to this basis is a matrix
we will denote by $S_8$, meaning that it is the $8 \times 8$ matrix
which has zeroes everywhere except for a line of ones connecting the
$(1,8)$ and the $(8,1)$ entries.  Also, the canonical involution $\bar{\ }$ is
given by
\[
\overline{u_i} = \begin{cases}
-u_i & \text{if $i \ne 4, 5$} \\
u_5 & \text{if $i = 4$} \\
u_4 & \text{if $i = 5$}.
\end{cases}
\]
\end{borel*}

\begin{borel}{Over $K$, $H$ is isomorphic to $\Spin_{10} \semid \Z/2$} \label{Adef}
Let $A$ denote the 10-dimensional subspace $e_0 \times J$ of $J$,
which is $A = \sjordmat{0}{F}{F}{\C}{0}{0}$.  We first observe that if
$f \in H$, then $f(e_0 \times j) = f^\dag(e_0) \times f^\dag(j)$, so
$f(A) = A$.  The multiplication on $J$ restricts to give $A$ the
structure of a central simple Jordan algebra as well, albeit with a
different unit element.  It has norm form $N_A$ given by
\[
N_A\sjordmat{0}{\alpha}{\beta}{c}{0}{0} = \alpha \beta - \n(c).
\]

Next we extend scalars to $F(t)$ and consider $N(t e_0 + j) = N(f(t
e_0 + j)) = N(t e_0 + f(j))$.  The coefficient of $t$ in this
expression is $T(e_0, j^\#) = T(e_0, f(j)^\#)$.  For $j$ actually
lying in $A$, $T(e_0, j^\#) = N_A(j)$, so $f$ must restrict to
preserve the norm on $A$.  We write $O(A)$ for the algebraic
subgroup of $GL(A)$ consisting of maps which preserve the norm $N_A$,
and we have proven that restriction provides a map 
$H \ra O(A)$ which is defined over $K$.  (It may not be defined over
$F$ because of the strange $\iota$-action on $H$ inherited from $E^K_6$.)

Clearly, $O(A)$ is an orthogonal group for the 10-dimensional
quadratic form $\hyp \perp -\n$.  Moreover, the map $H \ra O(A)$
has kernel of order 2: Anything in $H$ which maps to the identity in
$O(A)$ fixes all of the diagonal idempotents, and so must be of the
form $g_\ut$ for some $\ut \in \Spin(\n)$ \cite[p.~155, Thm.~1]{Soda}.
However, $t_0$  
must also be 
the identity, so $\ut = (1,1,1)$ or $(1, -1, -1)$ by \cite[1.5(2)]{G:iso}.

We would like to show that the map $H \ra O(A)$ is surjective.
Note that $O(A)$ is generated by 
\begin{itemize}
\item the special orthogonal group $SO(B)$ for $B$ the codimension 1
subspace of $A$ spanned by $\C$ and $e_1 - e_2$ (where $e_i$ denotes
the element of $J$ whose only nonzero entry is a 1 in the $(i+1,
i+1)$-position) endowed with the quadratic form given by restricting $N_A$; 
\item anything in $O(A)$ with determinant $-1$; and 
\item anything in $O(A)$ which doesn't leave $B$ invariant.
\end{itemize}
Since for $f \in \aut(J)$, $f^\dag = f$, the subgroup
$\aut(J/e_0)$ of elements of $\aut(J)$ which fix $e_0$ is a
subgroup of $H$.  As described in \cite[p.~376, Thm.~4]{Jac:J}, 
$\aut(J/ e_0) \cong \Spin(B)$ and the restriction to $B$ gives the
surjection onto $SO(B)$.  The map 
\[
\sbasjord \mapsto \sjordmat{\e_0}{\e_2}{\e_1}{c_0}{c_2}{c_1}
\]
lies in $H(K)$ and restricts to have determinant $-1$ on $A$.
Finally, we consider Freudenthal's maps from \cite[p.~74]{Jac:J3}.
For $E_{ij} \in M_3(\C)$ the matrix whose only nonzero entry is a 1
in the $(i,j)$-position, $1_3$ is the $3 \times 3$ identity matrix, $x
\in \C$, and $a \in J$, he defines a 
map $\psi_{ij}(x) \in \Inv(J)$ given by
\[
\psi_{ij}(x)(a) = (1_3 + xE_{ij}) a (1_3 + x E_{ij})^\ast,
\]
where juxtaposition denotes the usual multiplication in $M_3(\C)$, not
the Jordan multiplication.
So $\psi_{ij}(x) \in H(K)$ if $i, j \ne 1$.  In particular,
$\psi_{32}(u_5) \vert_A$ is given by
\[
\psi_{32}(u_5) \vert_A \stbtmat{\e_1}{c_0}{\cdot}{\e_2} = 
\stbtmat{\e_1}{c_0 - \e_1 u_4}{\cdot}{\e_2 - \n(c_0, u_4)},
\]
which doesn't leave $B$ invariant.

Finally, we observe that $\cc{H}$ is isomorphic to $\Spin(A)$.  The
inverse image, call it $H'$, of $SO(A)$ maps onto $SO(A)$ with a
kernel which is central and of order 2.  
Consequently, $H'$ is absolutely almost simple and
hence must be isomorphic to $\Spin(A)$.  Since $H'$ is connected and
$[G : H'] = 2$, $\cc{H} = H'$.
\end{borel}

\begin{borel}{Over $F$, $H$ is isomorphic to $\Spin(4\hyp \perp
\qform{-2, 2k}) \semid \Z/2$} \label{twistA}
To compute the isomorphism class of $H$ over $F$, we observe that the
map $h \mapsto h^\dag$ restricts to the identity on the kernel of the
$K$-map $H \ra O(A)$, so the $\iota$-action on $H$
induces one on $O(A)$, which we will calculate explicitly.

Fix the basis $(u_1, u_2, u_3, u_4, e_1, e_2, u_5, \ldots, u_8)$ for
$A$ so that the Gram matrix for $N_A$ becomes
\[
\left( \begin{smallmatrix}
& & -S_4 \\ 
& S_2 & \\
-S_4 & & 
\end{smallmatrix} \right),
\]
for $S_2$ and $S_4$ defined analogously to how $S_8$ was in \ref{udef}.
Then $SO(A)$ is generated by 
\begin{itemize}
\item a torus $T$ consisting of diagonal
matrices with diagonal entries $(d_1, d_2, \ldots, d_5, d_5^{-1},
d_4^{-1}, \ldots, d_1^{-1})$;

\item root groups $U_{ij} \!: \Ga \ra SO(A)$ given by 
\[
U_{ij}(r) = 1_{10} + rE_{ij} - rE_{j^\ast i^\ast}
\]
for $1_{10}$ the $10\times 10$ identity matrix, $i^\ast := 11 - i$, and
$(i, j) = (i, i+1)$ for $i = 1, 2, 3$, and  
their transposes; and

\item root groups $V_{ij} \!: \Ga \ra SO(A)$ given by 
\[
V_{ij}(r) = 1_{10} + r(E_{ij} + E_{j^\ast i^\ast})
\]
for $(i, j) = (4, 5)$ and $(4, 6)$, and their transposes.  (Note that
$V_{45}(r) = \psi_{32}(r u_5)\vert_A$ and $V_{46}(r) =
\psi_{23}(ru_4)\vert_A$ for $r \in F = \Ga(F)$.) 
\end{itemize}

Since the torus lies in the image of $\Rel\Cn$ and 
$g_\ut^\dag = g_{\s_\n(\ut)^{-1}}$, the action on $T$ and on the first kind
of root groups is the usual $\iota$-action.  However,
\[
V_{45}(r)^\dag = \psi_{32}(r u_5)^\dag = \psi_{23}(-r u_4) =
V_{46}(-r).
\]
So the map $h \mapsto h^\dag$ induces on $SO(A)$ the map $f \mapsto
MfM^{-1}$ for  
\[
M = \left( \begin{smallmatrix}
1_4 & & \\
& -S_2 & \\
& & 1_4
\end{smallmatrix} \right).
\]

Write $\eta$ for the 1-cocycle in $Z^1(K/F, O(A))$ given by
$\eta_\iota = M$.  The $K$-map $H \ra O(A)$ descends to a map over $F$
from $H$ onto the twisted group $O(A)_\eta$, so we wish to describe the
group $O(A)_\eta$.

But this is now just a problem of explicitly computing a quadratic
form given by descending down a quadratic extension.  So we need to
find a $K$-basis of $A \otimes K$ consisting of elements fixed by the
map $a \otimes \kappa \mapsto M(a) \otimes \iota(\kappa)$.  Then
$O(A)_\eta \cong O(q)$, where $q$ is the restriction of $N_A$ to the
$F$-span of those fixed vectors.  Such a $K$-basis is given by $u_i$
for $1 \le i \le 8$, 
$e_1 - e_2$, and $\sqrt{k} e_1 + \sqrt{k} e_2$.  These vectors
give an orthogonal basis for a quadratic form $4 \hyp \perp \qform{-2,
2k}$, which proves the claim.
\end{borel}

We close this section by proving a little lemma which
foreshadows the way we will prove the Main Theorem for quasi-split groups
of type $\dEvi$.

\begin{lem} \label{Hres}
The restriction of the Rost invariant on $H^1(F, E^K_6)$ to the image
of $H^1(F, \cc{H})$ has trivial kernel.
\end{lem}

\begin{pf}
We first observe that the embedding of $\cc{H}$ in $E^K_6$ has Rost
multiplier one.  Since the Rost multiplier is invariant under
scalar extension, we may work over $K$, where this embedding is
described in \ref{Adef}.  Then some of the coroots (identified with
copies of $\Gm$ lying in the maximal torus $T$ from \ref{twistA})
for $\cc{H}$ are the same as those for $\Spin(\n)$ considered as a
subgroup of $\Inv(J)$ via the map $g$.  Since the embedding $\Spin(\n)
\injects \Inv(J)$ has Rost multiplier one, so does $\cc{H} \injects
E^K_6$.

We are left with showing that anything in the kernel of the Rost invariant on
$H^1(F, \cc{H})$ maps to be trivial in $H^1(F, E^K_6)$.  By Lemma
\ref{structlem}, the kernel of the Rost invariant is precisely the
image of $H^1(F, C)$, where $C \cong \mmu2$ has for its unique
nontrivial element the map $g_{(1, -1, -1)}$ for $g$ as in \eqref{gdef}.
Now $\Spin(\n)$ is a subgroup of $\cc{H}$ over $F$ and it contains
this subgroup 
$C$ so we have a composition
\begin{neqn} \label{bigcomp}
H^1(F, C) \ra H^1(F, \Spin(\n)) \ra H^1(F, \cc{H}) \ra H^1(F, E^K_6).
\end{neqn}
Since $\n$ is hyperbolic, the spinor norm $SO(\n)(F) \ra \Fsq = H^1(F,
C)$ is surjective, so by the exactness of \eqref{Spin.coseq} the first
map in \eqref{bigcomp} is trivial.
\end{pf}

\section{Special cocycles} \label{specsec}

\begin{defn} \label{zdef}
For $\ua = \trion{a} \in \trio{(\Fx)}$ with product 1, we define a
``special'' cocycle $z := z_{K,\ua}$ in $H^1(K/F, G)$.  Set $z_\iota
= (z_0, z_1, z_2)$ where $z_j = m_j(\ua) d P$ for $P$ the permutation matrix
giving the map $u_k \mapsto u_{\pi(k)}$ for $\pi$ the permutation
$(1\,2)(3\,6)(4\,5)(7\,8)$, $m_j(\ua)$ the diagonal matrix 
\begin{neqn} \label{mjdef}
m_j(\ua) := \mathrm{diag}(1, a_j, a_j, a_{j+2}^{-1}, a_{j+1}^{-1},
1, 1, a_j)
\end{neqn}
with subscripts taken modulo 3, and
\[
d := \mathrm{diag}(1, 1, -1, 1, 1, -1, 1, 1).
\]
Then the $z_j$ form a related triple by
\cite[1.6, 1.7, 1.5(3)]{G:iso}, so $z_\iota \in G(K)$.  Since $P$ is
an isometry of $\n$,   
$\s_\n(P) = P^{-1} = P$.  Then 
\[
^\iota z_j = \s_\n(m_j(\ua) d P)^{-1} =  \diag(a_j^{-1}, 1, 1, a_{j+1},
a_{j+2}, a_j^{-1}, a_j^{-1}, 1) \, d \, P
\]
and so $z$ is indeed in $Z^1(K/F, G)$.
\end{defn}

\begin{borel}{Freedom in the definition} \label{freedom}
Of course, some of these special cocycles are cohomologically
equivalent in $H^1(K/F, G)$.  If $\ua$ and $\underline{a'}$ are two
triples in $\trio{(\Fx)}$ such that $a_j^{-1} a'_j \in
N_{K/F}(K^\ast)$ for all $j$, fix  
$\la_j \in K^\ast$ such that $a_j^{-1} a'_j = \la_j \iota(\la_j)$.
Then for $\ul = \trion{\ell}$ with $\ell_j =
P m_j(\underline{\lambda}) P$, $\ul$ is a related triple by
\cite{G:iso}, so $\ul \in G(K)$.  Then
$^\iota\hskip-.1em \ul \,\, (z_{K,\underline{a'}})_\iota \,\, \ul^{-1} = (z_{K,
\ua})_\iota$,
i.e., the two cocycles $z_{K,\ua}$ and $z_{K,\underline{a'}}$ are cohomologous.
\end{borel}

\begin{borel*} \label{sstwist}
We will twist by these cocycles to move a cocycle in $H^1(F, G)$ so that it takes values in a
semisimple group.  For now, we just observe that the semisimple group
we get from one of them, $\Spin(\n)_z$, is described in \cite[pp.~403,
404]{G:iso}:  
Let $k \in \Fx$ be such that $K = F(\sqrt{k})$ and let $Q_i$ 
denote the quaternion algebra $(k, a_i)_F$ generated by elements $x, y$ such
that $x^2 = k$, $y^2 = a_i$, and $xy = -yx$.  Then the group $\Spin(\n)_z$ is 
isomorphic to $\Spin (A_i, \s_i)$ where $A_i$ is isomorphic to
$M_4(Q_i)$, 
$\s_i$ is an isotropic orthogonal involution with trivial
discriminant, and
\begin{neqn} \label{tri.rel}
(C_0(A_i, \s_i),\underline{\s_i}) \cong (A_{i+1}, \s_{i+1}) \times
(A_{i+2}, \s_{i+2}),
\end{neqn}
where the subscripts are taken modulo 3.
(These properties specify the $\s_i$ up to isomorphism \cite[2.3]{G:clif}.)
\end{borel*}

The Moving Lemma lets us say something useful about the Rost
invariant of our special cocycles.

\begin{cor} \label{specialcor}
The Rost invariant of $z_{K,\ua}$ is trivial if and only if
the 1-cocycle is cohomologically trivial in $H^1(F, E^K_6)$.
\end{cor}

\begin{pf}
Consider the element $j = \sjordmat{0}{0}{0}{c}{0}{0}$ in $e_0 \times
J_K$ for $c = (u_2 + u_8) / 2$.  Then $\n(c) = 0$ and so consulting the
explicit formula for $j^\#$ in \cite[p.~358]{Jac:J}, we see that $j^\# = 0$.
Moreover, for $z := z_{K,\ua}$, we have $z_\iota \iota j =
\sjordmat{0}{0}{0}{c'}{0}{0}$ for $c' = (u_1 + u_7) / 2$.  Then 
\[
T(j, j') = c \overline{c'} + \overline{c} c' = 2 \n(c, c') = 1.
\]
Applying the Moving Lemma shows that $z$ is equivalent in
$H^1(K/F, E^K_6)$ to some $z' \in Z^1(K/F, G)$ whose image in
$H^1(K/F, P)$ is $(1, a_0^{-1}, a_0)$.  So in particular we may assume
that $\s_\n(z'_\iota) z'_\iota = (1, a_0^{-1}, a_0)$.
Since the restriction of $z'_\iota$
to the 10-dimensional subalgebra $A$ defined in \ref{Adef} has
determinant 1, it must lie in $\cc{H}$.
The statement about the triviality is now a consequence of Lemma \ref{Hres}
\end{pf}

In a special case we can calculate the value
of the Rost invariant of our special cocycles quite explicitly.

\begin{lem} \label{rostcalc}
For $a, k \in \Fx$ such that $K = F(\sqrt{k})$,
the Rost invariant of the 1-cocycle $z_{K,(1, a, a^{-1})}$ is $(a)
\cup (k) \cup (-1)$.
\end{lem}

\begin{pf}
The cocycle $z := z_{K, (1, a, a^{-1})}$ takes values in $H$ and
 restricts to have determinant one on the subalgebra $A$ defined in
 \ref{Adef}, so $z \in Z^1(K/F, \cc{H})$.  
Since the embedding $\cc{H} \subset E^K_6$ has Rost
multiplier 1, to compute the Rost invariant of $z$, we may compute the
Rost invariant of $z$ in $H^1(F, \cc{H})$.  But recall that $\cc{H} \cong
 \Spin(q)$ for $q = 4\hyp \perp \qform{-2, 2k}$ and that $H^1(F,
 SO(q))$ classifies nondegenerate quadratic forms of the same
 dimension and discriminant as $q$.  So we can compute the Rost
 invariant of $z$ by computing the quadratic form $q_z$ corresponding
 to the image of $z$ in $H^1(F, SO(q))$, which is just the restriction
 of $q \otimes K$ to the vector subspace fixed by the action $a
 \otimes \kappa \mapsto z_\iota M(a) \otimes \iota(\kappa)$ for $M$ as
 in \ref{twistA}.

We will perform the Galois descent calculation by decomposing $A
\otimes K$ into 2-dimensional subspaces and calculating the Galois action
on those subspaces.
\[
\begin{array}{cccc}
\text{subspace} & \text{restriction of} &  \text{$F$-basis for}       &
\text{contribution} \\
\text{basis}    & \text{$z_\iota M$} & \text{fixed subspace} & \text{to
$q_z$} \\ \hline
(u_1, u_2) & S_2  & & \text{totally}\\
(u_7, u_8) & S_2  & & \text{isotropic}\\
(u_3, u_6) & -S_2  & 
u_3 - u_6, \sqrt{k} u_3 + \sqrt{k} u_6 & \qform{2, -2k}\\
(u_4, u_5) & \stbtmat{}{a}{a^{-1}}{} & 
a u_4 + u_5, -a\sqrt{k} u_4 + \sqrt{k}u_5 & \qform{-2a, 2ak} \\
(e_1, e_2) & \stbtmat{}{-a^{-1}}{-a}{}  & 
-e_1 + a e_2, \sqrt{k} e_1 + a \sqrt{k} e_2 & \qform{-2a, 2ak} 
\end{array}
\]
The first two subspaces form a complementary pair of totally isotropic
subspaces, so they contribute two hyperbolic planes to $q_z$.
Thus the image of $z$ is $q_z = 2\hyp \perp \qform{2, -2k, -2a, 2ak,
-2a, 2ak}$ and the Rost invariant of $z$ is the Arason invariant of 
$q_z - q = \qform{2}\pform{a, k, -1}.$
\end{pf}
\section{Quasi-split groups of type $E_6$ and $E_7$} \label{E6sec}

This section consists solely of a proof of the main theorem, beginning
with a nearly trivial lemma.

\begin{lem}\label{homlem}
Suppose that $C$ is a central subgroup in an almost simple simply
connected group 
$\Gamma$.  Then $H^1(F, C)$ acts on $H^1(F, \Gamma)$ and for $\zeta \in H^1(F,
C)$ and $\gamma \in H^1(F, \Gamma)$, we have
\[
R_\Gamma(\zeta \cdot \gamma) = R_\Gamma(\zeta) + R_\Gamma(\gamma),
\]
where $R_\Gamma(\zeta)$ denotes the image of $\zeta$ under the composition
$H^1(F, C) \longto H^1(F, \Gamma) \xrightarrow{R_\Gamma} H^3(F, \QZt)$.
\end{lem}

\begin{pf}
We pick a 1-cocycle $z \in Z^1(F, C)$ which represents $\zeta$.  Then
we have a diagram
\[
\begin{CD}
H^1(F, \Gamma) @= H^1(F, \Gamma_z) @>{\sim}>{\tau_z}> H^1(F, \G) \\
@V{R_\G}VV @V{R_{\G_z}}VV @VV{R_\G}V \\
H^3(F, \QZt) @= H^3(F, \QZt) @>>> H^3(F, \QZt) 
\end{CD}
\]
where the group $\G_z$ is the usual twist of $\G$ by the cocycle $z$, so
that it is just the group $\G$ with a different Galois action
so that a member $\s$ of the Galois group maps $g \mapsto z_\sigma \,
{^\sigma g} z^{-1}_\s$.  In our case, $z_\s$ is always central, so in
fact $\G_z = \G$.  The map $\tau_z$ is the usual twisting map
\cite[I.5.5]{SeCG}, defined by sending $a \in Z^1(F, \G_z)$ to the
1-cocycle $\s \mapsto a_\s z_\s$.  The composition of the two maps on
the top row is then the action of $\zeta$.

The left-hand box commutes because the Rost invariant is canonical.
The right-hand box commutes, where the bottom map is given by $\eta
\mapsto \eta + R_\G(\zeta)$ by \cite[p.~16, Lem.~8]{Gille:inv}.
\end{pf}

This result has the obvious corollary that
the induced map
$H^1(F, C) \longto H^3(F, \QZt)$ is a group homomorphism.

\begin{borel}{Groups of type $\oEvi$}
Suppose first that our simply connected quasi-split group of type
$E_6$ is split and denote it simply by $E_6$.  From Example
\ref{e6sm}, we have an embedding $F_4 \times \mmu3 \injects E_6$ which
induces a surjection on $H^1$ terms.  So for $\varepsilon \in H^1(F,
E_6)$, we can find a $\phi \in H^1(F, F_4)$ and $\zeta \in H^1(F,
\mmu3)$ such that $\phi \oplus \zeta \mapsto \varepsilon$.  Now we
simply note that
since $E_6$ is 
split, the image of $H^1(F, \mmu3) \ra H^1(F, E_6)$ is trivial.  So if
$\varepsilon$ is in the kernel of the Rost invariant $R_{E_6}$, by
Lemma \ref{homlem} $\phi$ must be killed by the composition
\[
H^1(F, F_4) \ra H^1(F, E_6) \xrightarrow{R_{E_6}} H^3(F, \QZt).
\]
As described in \ref{e6rm}, the Rost multiplier of the embedding $F_4
\subset E_6$ is 1, so $\phi$ lies in the kernel of the Rost invariant
$R_{F_4}$, which is known to be trivial.  So $\varepsilon$ is the
image of $\zeta$, which we have already observed is trivial.
\end{borel}

\begin{rmk}[Noninjectivity for $\oEvi$]
We want to provide an example to show that the Rost invariant is
noninjective for the group $E_6$.  We can not simply apply Remark
\ref{noninj.F4} and the fact that the embedding $F_4 \injects E_6$ has
Rost multiplier one, since two isotopic Albert algebras have the same
image in $H^1(F, E_6)$.

So instead fix a ground field $F$ which supports a division (=
nonreduced) Albert $F$-algebra $J$.  Then over the field $F(t)$, the
norm $N$ of $J$ doesn't represent $t$ as can be seen by some elementary
valuation theory \cite[p.~417, Lem.~1]{Jac:J}.  Consequently, $N$ is
not isomorphic to $tN$ over $F(t)$, so the images of the
two classes $(J) \oplus (1)$ and $(J) \oplus (t)$ under the map
$H^1(F, F_4) \times H^1(F, \mmu3) \ra H^1(F, E_6)$ are distinct by
\cite[2.8(2)]{G:struct}.  However, since the image of $H^1(F, \mmu3)
\ra H^1(F, E_6)$ is trivial, by Lemma \ref{homlem} the two classes in
$H^1(F, E_6)$ have the
same Rost invariant.
\end{rmk}

\begin{borel}{Groups of type $\dEvi$}
Suppose now that our quasi-split simply connected group of type $E_6$
is not actually split, so that it only becomes split over some
quadratic field extension $K$ of $F$.  We write $E^K_6$ for this
group, as we have since Section \ref{Gsec}.  
Then by the split case, any $\alpha \in H^1(F, E^K_6)$ which
is in the kernel of the Rost invariant must become trivial over $K$
and so must come from $H^1(K/F, E^K_6)$.  Applying Ferrar's Lemma
\ref{ferrlem}, we have that $\alpha$ is the image of some $\beta \in
H^1(K/F, G)$.  
\end{borel}

\begin{borel}{Twisting} \label{clevertwist}
We fix a triple $\ua = \trion{a} \in \trio{(\Fx)}$ such that $a_0 a_1
a_2 = 1$ which represents the 
image of $\beta$ in $H^1(K/F, P)$.  (This makes sense thanks to the
description of $H^1(K/F, P)$ 
in \ref{Pdesc}.) 
Then we set $z := z_{K,\ua}$ as defined in \ref{zdef},
and we can twist $E^K_6$ by $z$ to obtain a diagram 
\[
\begin{CD}
H^1(F, G_z) @>>> H^1(F, (E^K_6)_z) @>{R_{(E^K_6)_z}}>> H^3(F, \QZt) \\
@V{\tau_z}VV @V{\tau_z}VV @VV{\cdot + R_{E^K_6}(z)}V \\
H^1(F, G) @>>> H^1(F, E^K_6) @>{R_{E^K_6}}>> H^3(F, \QZt),
\end{CD}
\]
where the right vertical arrow has the specified value by
\cite[p.~16, Lem.~8]{Gille:inv}. 
\end{borel}

\begin{borel}{The image of $\tau_z^{-1}(\beta)$ in $H^1(F, SO\As)$}
\label{SOim} 
We want to say something about what kind of class $\beta' :=
\tau_z^{-1}(\beta)$ can be.  In particular, its image in 
$H^1(K/F, P_z)$ is trivial, so $\beta'$ comes from the
semisimple part of $G_z$, which is isomorphic to $\Spin\As$ for $\As$
one of the three algebras $A_i$ described in \ref{sstwist}.   

So we can think of $\beta'$ as lying in $H^1(K/F, \Spin\As)$ and consider 
its image in $H^1(K/F, SO\As)$.  Let $L$ be a generic splitting field
of $A$ (e.g., a function field of its Severi-Brauer variety) and
consider the image of $\beta'$ in $H^1(L, SO\As)$.  Since $A$ is split
by $L$, $\s$ becomes adjoint to the quadratic form $\pform{k, a_{i+1}}
\perp 2\hyp$ \cite[2.3]{G:clif}.  The image of $\beta'$ determines an
8-dimensional 
quadratic form $q$ over $L$, and the Rost invariant of $\beta'$ is
just the class of $q - \pform{k, a_{i+1}}$ in $\IL3 / \IL4$.
However, by the twisting argument above, the Rost invariant
of $\beta'$ over $F$ is $-R_{E^K_6}(z)$.  Since $A$ is split over $L$,
$a_i \in \Lx$ is a norm from $KL$, so by \ref{freedom} and   
Lemma \ref{rostcalc} the Rost invariant becomes $(k) \cup (a_{i+1})
\cup (-1)$ over $L$.  

So we have that for $\phi = \pform{k, a_{i+1}}$, $q - \phi \in \IL3$ and
$q - \phi \equiv \phi \pform{-1} \mod{\IL4}$.  But then 
\[
q + \phi = (q - \phi) + 2\phi \equiv  4 \phi \equiv 0 \mod \IL4.
\]
So $q + \phi \in \IL4$.  However $\dim(q \perp \phi) = 12 < 16$, so by
the Arason-Pfister Hauptsatz, $q \perp \phi$ is hyperbolic and $q
\cong \qform{-1} \phi \perp 2\hyp$.

The consequence of all that is that the image of $\beta'$ in $H^1(L,
SO\As)$ is the same as the image of $-1 \in \Fsq = H^1(F, Z(SO\As))$.
Since $A$ is Brauer-equivalent to a quaternion algebra, it follows
from the material in \cite[Ch.~10]{Sch} that
the canonical map $H^1(F, SO\As) \ra H^1(L,
SO\As)$ is injective.  (This was shown independently in \cite{PSS:herm} and
\cite[5.3.1]{Dej:th}.)
So the image of $\beta'$ in $H^1(F, SO\As)$ must also be $-1$.  
\end{borel}

\begin{borel*}
More generally, any simply connected group $\Gamma$ of type $\oD$ is
isomorphic to $\Spin(A_i, \s_i)$ for three central simple algebras
$A_i$ of degree 8 with $i = 0, 1, 2$ endowed with an orthogonal
involution $\s_i$ with trivial discriminant and related as in
\eqref{tri.rel}.  

Each of the three descriptions of $\Gamma$ comes
paired with natural maps $\Gamma \ra SO\Asi \ra P\Gamma$ for $P\Gamma$
the adjoint group associated to $\Gamma$.  The kernel of the second
map is $Z(SO\Asi) \cong \mmu2$, and the kernel of the composition is
$Z(\Gamma)$, which is isomorphic to the subgroup of $\mmut2$ with
product one.  Then the group $H^1(F, Z(\Gamma))$ can be identified
with the set of triples $\ub = \trion{b} \in \Fsq$ with product 1
\cite[44.14]{KMRT} and where the map $H^1(F, Z(\Gamma)) \ra H^1(F,
Z(SO\Asi))$ is given by $\ub \mapsto b_i$.

\begin{lem} (Notation as in the preceding paragraph.) \label{trilem}
Suppose an element $\eta \in H^1(F, \Gamma)$ has the same image in
$H^1(F, SO\Asi)$ as $c_i \in \Fsq = H^1(F, Z(SO\Asi))$ for $i = 1, 2$.  
Then $\eta$ is the image of $((c_1 c_2)^{-1}, c_1, c_2)$ coming
from $H^1(F, Z(\Gamma))$.
\end{lem}

\begin{pf}
We have a short exact sequence $1 \ra Z(SO\Asi) \ra SO\Asi \ra P\Gamma
\ra 1$, so $\eta$ is killed by the composition $H^1(F, \Gamma) \ra
H^1(F, SO\Asi) \ra H^1(F, P\Gamma)$ for $i = 1$. 
Thus $\eta$ lies in the image of $H^1(F, Z(\Gamma))$.

For general Galois-cohomological reasons, the map $H^1(F, Z(\Gamma))
\ra H^1(F, \Gamma)$ is a group homomorphism.  (Although the second set
doesn't have a group structure, the image of the first set does.)
The kernel of this map 
can be described fully by suitably applying \cite[35.4]{KMRT},
but for our purposes it is enough to observe that it contains all
elements of the form $(s, s^{-1}, 1)$ for $s$ a spinor norm of an element in
$SO(A_2, \s_2)(F)$ and symmetrically.  Let $\cc{G\Asi}$ be the
algebraic group of {\em proper similarity factors},
i.e., the group with $F$-points
\[
\cc{G\Asi}(F) = \left\{ m \in \Fx \mid 
\text{$\exists \, f \in A_i^\ast$ such that $m = \s_i(f)f$ and
$\Nrd_{A_i}(f)= m^4$}
\right\}.
\]
Then for every $m_0 \in \cc{G(A_0, \s_0)}(F)$, the kernel contains an
element of the form $\trion{m}$ and symmetrically.  Conversely,
if $\trion{b}$ is in the kernel, then $b_i \in \cc{G(A_i, \s_i)}(F)$
for all $i$.  

It is also the case that the natural map $\Fsq = H^1(F, Z(SO\Asi)) \ra H^1(F,
SO\Asi)$ is a group homomorphism, and its kernel is precisely
$\cc{G\Asi}(F)$.  So the fact that $\eta$ maps to the image of $c_2$
in $H^1(F, SO(A_2, \s_2))$ 
means that $\eta$ is the image of some class $\trion{n}$ in $H^1(F,
Z(\Gamma))$ such that $n_2 = c_2$.

Now consider the middle component of this triple.  By hypothesis,
$n_1 = m_1 c_1$ for some $m_1 \in \cc{G(A_1, \s_1)}(F)$.
By \cite[p.~262, Prop.]{M:norm}, the group $\SN(A_2, \s_2)(F)$ of spinor norms 
from $SO(A_2, \s_2)(F)$ is $\Fxsq$ times the norms from field extensions $E$
which split $A_2$ and make $\s_2$ isotropic.  
By \cite[p.~263, Prop.]{M:norm}, $\cc{G(A_1, \s_1)}(F)$ is equal to the group
generated by the norms from every extension field $E$ which splits
$A_1$ and makes $\s_1$ hyperbolic.  
Since the $\Asi$ are related by \eqref{tri.rel}, any
extension which splits $A_1$ and makes $\s_1$ hyperbolic certainly
splits $A_2$ and makes $\s_2$ isotropic, so $\SN(A_2, \s_2)(F) 
\supseteq \cc{G(A_1, \s_1)}(F)$.  Consequently, $\trion{n}$ has the
same image in $H^1(F, \Gamma)$ as $((c_1 c_2)^{-1}, c_1, c_2)$.
\end{pf}
\end{borel*}

\begin{borel}{$\beta'$ is in the image of $H^1(K/F, Z(\Spin\As))$}
Let $\As = (A_0, \s_0)$ for $\Asi$ as in \ref{sstwist}.  Combining the
result from \ref{SOim} with Lemma \ref{trilem}, we have that $\beta' \in
H^1(F, \Spin\As)$ is the image of $(1, -1, -1) \in H^1(F,
Z(\Spin\As))$.  However, for $k \in \Fx$ such that $K = F(\sqrt{k})$,
since $K$ certainly splits $A$ and makes $\s$ hyperbolic and $-k =
N_{K/F}(\sqrt{k})$, by Merkurjev's norm principle \cite[p.~262,
Prop.]{M:norm} there is some 
element of $SO\As(F)$ with spinor norm $-k$.  Then as described in the
proof of Lemma \ref{trilem}, $\beta'$ is also the
image of $(1, k, k^{-1}) \in H^1(F, Z(\Spin\As))$, which itself is in
the image of $H^1(K/F, Z(\Spin\As))$.
\end{borel}

\begin{borel*} \label{closer.2E6}
Consider the 1-cocycle $b = \tau_z(b') \in Z^1(K/F, G)$ for $b'$ the
image of $(1, k, k^{-1})$ as above.  (Note that $b$ represents the
class of $\beta$ and is the 1-cocycle which takes the value $g_{(1,
-1, -1)} z_{K,\ua}$ at $\iota$.)  For $j$ and $c$ as in the proof of
\ref{specialcor}, we set $j' := 
b_\iota \iota j$, so that $j' = \sjordmat{0}{0}{0}{c'}{0}{0}$ for $c'
= (u_1 + u_7) / 2$.  
So $T(j, j') = 2\n(c,c') = 1$.  
Then by the Moving Lemma \ref{movlem}, we may
replace $\beta$ by a different inverse image of $\alpha$ in $H^1(K/F,
G)$ and so assume that $\ua = (1, a_0, a_0^{-1})$.

Any element of $G$ with multiplier $(1, \cdot, \cdot)$ lies in $H$,
and since such an element restricts to have determinant one on the
subspace $A$ defined in  
\ref{Adef}, it in fact lies in $\cc{H}$.  Thus $\alpha$ is in the
image of $H^1(F, \cc{H})$.  Since the Rost
invariant of $\alpha$ is trivial, $\alpha$ must be the trivial class
by Lemma \ref{Hres}.
\end{borel*}

\begin{borel}{Groups of type $E_7$}
We are left with proving that the Rost invariant has trivial kernel
for $G$ split of type $E_7$, but 
this follows directly from the same conclusion for quasi-split groups
of type $E_6$, 
thanks to Proposition
\ref{e7smprop}.
\end{borel}


\subsection*{Acknowledgements}

I thank 
ETH in Z\"urich for its hospitality while I performed part of the
research contained herein;
Ph.~Gille and Zinovy Reichstein for helpful discussions; Vladimir
Chernousov for catching an error in an earlier version of this paper; 
and absolutely most of all Markus Rost for
explaining so much about his invariants and for posing this problem to me.

\providecommand{\bysame}{\leavevmode\hbox to3em{\hrulefill}\thinspace}

\bigskip

%
%
\noindent R.~Skip Garibaldi\\%
e-mail: {\tt skip@member.ams.org} \\%
web: {\tt http://www.math.ucla.edu/\~{}skip/}\\%

\noindent UCLA\\%
Dept.~of Mathematics\\%
Los Angeles, CA 90095-1555


\begin{thebibliography}{KMRT98}

\bibitem[AH81]{AH}
B.N. Allison and W.~Hein, \emph{Isotopes of some nonassociative algebras with
  involution}, J. Algebra \textbf{69} (1981), no.~1, 120--142.

\bibitem[All94]{A:survey}
B.N. Allison, \emph{Structurable algebras and the construction of simple {L}ie
  algebras}, Jordan algebras (Berlin) (W.~Kaup, K.~McCrimmon, and H.P.
  Petersson, eds.), de Gruyter, 1994, (Proceedings of a conference at
  Oberwolfach, 1992).

\bibitem[Bor91]{Borel}
A.~Borel, \emph{Linear algebraic groups}, second ed., Graduate Texts in
  Mathematics, vol. 126, Springer-Verlag, New York, 1991.

\bibitem[Bou68]{Bou:g4}
N.~Bourbaki, \emph{{G}roupes et alg\`ebres de {L}ie. {C}h. {IV}, {V}, {VI}},
  Hermann, Paris, 1968, Actualit\'es Scientifiques et Industrielles, No. 1337.

\bibitem[Bro69]{Brown:E7}
R.B. Brown, \emph{Groups of type ${E}_7$}, J. Reine Angew. Math. \textbf{236}
  (1969), 79--102.

\bibitem[Che00]{Ch:rinvlet}
V.I. Chernousov, private letter, October 2000.

\bibitem[Dej99]{Dej:th}
I.~Dejaiffe, \emph{Somme orthogonale d'alg\`ebres centrales \`a involution:
  propri\'et\'es et applications}, DSc thesis, June 1999, Universit\'e
  Catholique de Louvain, Belgium.

\bibitem[{\`E}la72]{Elash}
A.G. {\`E}la{\v{s}}vili, \emph{Canonical form and stationary subalgebras of
  points in general position for simple linear {L}ie groups}, Functional Anal.
  Appl. \textbf{6} (1972), no.~1, 44--53, [Russian original: Funkcional. Anal.
  i Prilo\v zen. \textbf{6} (1972), no. 1, 51--62].

\bibitem[Fer69]{Ferr:E6}
J.C. Ferrar, \emph{Lie algebras of type ${E}_6$}, J. Algebra \textbf{13}
  (1969), 57--72.

\bibitem[Fer72]{Ferr:strict}
J.C. Ferrar, \emph{Strictly regular elements in {F}reudenthal triple systems},
  Trans. Amer. Math. Soc. \textbf{174} (1972), 313--331.

\bibitem[Fer80]{Ferr:class}
J.C. Ferrar, \emph{On the classification of {F}reudenthal triple systems and
  {L}ie algebras of type ${E}_7$}, J. Algebra \textbf{62} (1980), no.~2,
  276--282.

\bibitem[Gara]{G:clif}
R.S. Garibaldi, \emph{Clifford algebras of hyperbolic involutions}, to appear
  in Math.~Z.

\bibitem[Garb]{G:struct}
R.S. Garibaldi, \emph{Structurable algebras and groups of type ${E}_6$ and
  ${E}_7$}, to appear in J.~Algebra.

\bibitem[Gar98]{G:iso}
R.S. Garibaldi, \emph{Isotropic trialitarian algebraic groups}, J. Algebra
  \textbf{210} (1998), 385--418.

\bibitem[Gila]{Gille:sc}
Ph. Gille, \emph{Cohomologie galoisienne des groupes algebriques quasi-deployes
  sur des corps de dimension cohomologique $\le 2$}, $K$-Theory, to appear.

\bibitem[Gilb]{Gille:inv}
Ph. Gille, \emph{Invariants cohomologiques de {R}ost en caract\'eristique
  positive}, Compositio Math., to appear.

\bibitem[Gil00]{Gille:E8}
Ph. Gille, April 2000, private letter.

\bibitem[GV78]{GV}
V.~Gatti and E.~Viniberghi, \emph{Spinors of $13$-dimensional space}, Adv. in
  Math. \textbf{30} (1978), no.~2, 137--155.

\bibitem[Hec84]{Heck}
A.~Heck, \emph{Involutive automorphisms of root systems}, J. Math. Soc. Japan
  \textbf{36} (1984), no.~4, 643--658.

\bibitem[Igu70]{Igusa}
J.~Igusa, \emph{A classification of spinors up to dimension twelve}, Amer. J.
  Math. \textbf{92} (1970), 997--1028.

\bibitem[Jac59]{Jac:J1}
N.~Jacobson, \emph{Some groups of transformations defined by {J}ordan algebras.
  {I}}, J. Reine Angew. Math. \textbf{201} (1959), 178--195.

\bibitem[Jac61]{Jac:J3}
N.~Jacobson, \emph{Some groups of transformations defined by {J}ordan algebras.
  {III}}, J. Reine Angew. Math. \textbf{207} (1961), 61--85.

\bibitem[Jac68]{Jac:J}
N.~Jacobson, \emph{Structure and representations of {J}ordan algebras}, AMS
  Coll.\ Pub., vol.~39, AMS, Providence, RI, 1968.

\bibitem[Jac71]{Jac:ex}
N.~Jacobson, \emph{Exceptional {L}ie algebras}, Lecture notes in pure and
  applied mathematics, vol.~1, Marcel-Dekker, New York, 1971.

\bibitem[Kim88]{Kim:surv}
T.~Kimura, \emph{A classification theory of prehomogeneous vector spaces},
  Representations of Lie groups, Kyoto, Hiroshima, 1986, Academic Press,
  Boston, MA, 1988, pp.~223--256.

\bibitem[KMRT98]{KMRT}
M.-A. Knus, A.S. Merkurjev, M.~Rost, and J.-P. Tignol, \emph{The book of
  involutions}, Colloquium Publications, vol.~44, AMS, Providence, RI, 1998.

\bibitem[Lam73]{Lam}
T.-Y.\ Lam, \emph{The algebraic theory of quadratic forms}, Benjamin, Reading,
  MA, 1973.

\bibitem[McC69]{McC:FST}
K.~McCrimmon, \emph{The {F}reudenthal-{S}pringer-{T}its constructions of
  exceptional {J}ordan algebras}, Trans. Amer. Math. Soc. \textbf{139} (1969),
  495--510.

\bibitem[Mer96]{M:norm}
A.S. Merkurjev, \emph{A norm principle for algebraic groups}, St.\ Petersburg
  Math.\ J. \textbf{7} (1996), no.~2, 243--264.

\bibitem[MS91]{MS:norm3}
A.S. Merkurjev and A.A. Suslin, \emph{The norm residue homomorhpism of degree
  3}, Math. USSR-Izv. \textbf{36} (1991), 349--368.

\bibitem[Pop80]{Pop:14}
V.L. Popov, \emph{A classification of spinors of dimension fourteen}, Trans.
  Mosc. Math. Soc. (1980), no.~1, 181--232, [Russian original: Trudy Moskov.
  Mat. Obshch. {\bf 37} (1978), 173--217].

\bibitem[Pop94]{Pop:sec}
V.L. Popov, \emph{Sections in invariant theory}, The Sophus Lie Memorial
  Conference (Oslo, 1992), Scand. Univ. Press, Oslo, 1994, pp.~315--361.

\bibitem[PR94]{PlatRap}
V.P. Platonov and A.~Rapinchuk, \emph{Algebraic groups and number theory},
  Academic Press Inc., Boston, MA, 1994, Translated from the 1991 Russian
  original by Rachel Rowen.

\bibitem[PSS]{PSS:herm}
R.~Parimala, R.~Sridharan, and V.~Suresh, \emph{Hermitian analogue of a theorem
  of {S}pringer}, preprint.

\bibitem[Ros99]{Rost:14}
M.~Rost, \emph{On the {G}alois cohomology of ${\Spin(14)}$}, preprint, March
  1999.

\bibitem[Sch66]{Schfr}
R.~Schafer, \emph{An introduction to nonassociative algebras}, Pure and Applied
  Math., vol.~22, Academic Press, New York-London, 1966.

\bibitem[Sch69]{Schatt}
D.J. Schattschneider, \emph{On restricted roots of semi-simple algebraic
  groups}, J. Math. Soc. Japan \textbf{21} (1969), 94--115.

\bibitem[Sch85]{Sch}
W.~Scharlau, \emph{Quadratic and hermitian forms}, Grund.\ math.\ Wiss., vol.
  270, Springer, Berlin-New York, 1985.

\bibitem[Ser62]{Se62}
J.-P. Serre, \emph{Cohomologie galoisienne des groupes alg\'ebriques
  lin\'eaires}, Colloq. Th\'eorie des Groupes Alg\'ebriques (Bruxelles, 1962),
  Librairie Universitaire, Louvain, 1962, pp.~53--68.

\bibitem[Ser94]{SeCG}
J.-P. Serre, \emph{Cohomologie {G}aloisienne}, fifth ed., Lecture Notes in
  Mathematics, vol.~5, Springer-Verlag, Berlin, 1994, [English translation:
  {\em Galois cohomology}, Springer, 1997].

\bibitem[Ser95]{SeCG:p}
J.-P. Serre, \emph{Cohomologie {G}aloisienne: progr\`{e}s et probl\`{e}mes},
  Ast\'erisque (1995), no.~227, 229--257, S\'eminaire Bourbaki, vol.\ 1993/94,
  Exp.\ 783.

\bibitem[Sod66]{Soda}
D.~Soda, \emph{Some groups of type ${D}_4$ defined by {J}ordan algebras}, J.
  Reine Angew. Math. \textbf{223} (1966), 150--163.

\bibitem[Spr62]{Sp:cubic}
T.A. Springer, \emph{Characterization of a class of cubic forms}, Nederl.\
  Akad.\ Wetensch. \textbf{65} (1962), 259--265.

\bibitem[SV68]{SpV}
T.A. Springer and F.D. Veldkamp, \emph{On {H}jelmslev-{M}oufang planes}, Math.
  Z. \textbf{107} (1968), 249--263.

\end{thebibliography}
\end{document}